\documentclass{ws-m3as}
\usepackage{epsfig}
\usepackage{colortbl}
\usepackage{ulem}

\def\@xthm#1#2{\@beginassumption{#2}{\csname the#1\endcsname}{}\ignorespaces}
\def\@ythm#1#2[#3]{\@opargbeginassumption{#2}{\csname the#1\endcsname}{#3}\ignorespaces}%
\def\@beginassumption#1#2#3{\par\addvspace{8pt plus3pt minus2pt}%
              \noindent{\csname#1headfont\endcsname#1\ \ignorespaces#3 #2.}%
              \csname#1font\endcsname\hskip.5em\ignorespaces}
\def\@endassumption{\par\addvspace{8pt plus3pt minus2pt}\@endparenv}
\usepackage{tikz}
\usetikzlibrary{shapes,arrows}

\usepackage{multirow}

\usepackage{hyperref}

\usetikzlibrary{shapes,arrows,decorations.markings}
\usetikzlibrary{decorations.pathmorphing}

\tikzstyle{block_long} = [rectangle, draw, fill=blue!20,
    text width=10.0em, text centered, rounded corners, minimum height=3em]
\tikzstyle{block_medium} = [rectangle, draw, fill=blue!20,
    text width=6.0em, text centered, rounded corners, minimum height=3em]
\tikzstyle{block} = [rectangle, draw, fill=blue!20,
    text width=3.0em, text centered, rounded corners, minimum height=3em]
\tikzstyle{line} = [thick, draw, dashed,  -stealth']

\usepackage{mathrsfs,bbm}
\usepackage[american]{babel}

\newcommand{\R}{\mathbb{R}}

\newcommand{\cA}{\mathcal{A}}

\newcommand{\cD}{\mathcal{D}}

\newcommand{\cP}{\mathcal{P}}

\newcommand{\cM}{\mathcal{M}}

\newcommand{\EE}{\mathbb{E}}

\newcommand{\vf}{\varphi}

\newcommand{\RR}{\mathbb{R}}

\newcommand{\bomega}{\boldsymbol{\omega}}

\newcommand{\bu}{{\boldsymbol{u}}}
\newcommand{\bx}{{\boldsymbol{x}}}

\newcommand{\bv}{{\boldsymbol{v}}}

\newcommand{\bF}{\boldsymbol{F}}

\newcommand{\bbf}{\boldsymbol{f}}

\def\va{\raise 2pt\hbox{,}}

\def\cG{{\cal G}}

\def\cL{{\cal L}}

\def\cA{{\cal A}}

\def\cD{{\cal D}}

\newcommand{\p}{\partial}

\input epsf
\begin{document}


\markboth{N.~Bellomo, D.~Burini and J. Liao}{On the Way to the Complexity of Living Systems}

%
\catchline{}{}{}{}{}
%

\title{New Trends in Kinetic Theory \\ Towards the Complexity of Living Systems}


%

\author{Nicola Bellomo$^{(1),(2)}$, Diletta Burini$^{(3)}$, and Jie Liao$^{(4)}$}

\address{$^{(1)}$Departamento de Matem\'atica Aplicada and  Research Unit ``Modeling Nature'' (MNat), Facultad de Ciencias, Universidad de Granada, 18071 Granada, Spain.\\
nicola.bellomo@ugr.es}

\address{$^{(2)}$Polytechnic University of Torino, Italy.\\
nicola.bellomo@polito.it}

\address{$^{(3)}$Department of Mathematics and Computer Sciences, University of Perugia, Italy.\\
diletta.burini@unipg.it}

\address{$^{(4)}$School of Mathematics, Shanghai University of Finance and Economics,\\ Shanghai 200433, China.  \\
liaojie@mail.shufe.edu.cn}

\vskip1truecm

\maketitle

\begin{history}
\received{(Day Month Year)} \revised{(Day Month Year)}
\comby{(xxxxxxxxxx)}
\end{history}

\begin{abstract}
The development of a mathematics for living systems is one of the most challenging prospects of this century. The search began with the pioneering contribution of Ilia Prigogine, who developed methods from statistical physics to describe the dynamics of vehicular traffic. This visionary seminal research contribution has given rise to a great deal of research activity, which began at the end of the last century and has been further developed in this century by several authors who have developed mathematical methods, generally focused on applications. These methods are somewhat inspired by the classical kinetic theory, but significant differences have led to the concept of active particles and to a kinetic theory that is ultimately very different from the classical theory. Different approaches have been developed, each of which is in some way an alternative to the others. This paper develops a critical analysis of the scientific activity after Prigogine with the aim of developing a unified mathematical theory, taking into account the conceivable interactions that a mathematical theory of living systems can have with studies of artificial intelligence.
\end{abstract}



\renewcommand{\arraystretch}{1.2}

\parindent=20pt

\newpage

 \tableofcontents

\section{Motivations and plan of the paper}\label{Sec:1}

Mathematical methods inspired by the kinetic theory of classical particles, were introduced by Ilia Prigogine in his celebrated book, with Hermann, on the modeling of vehicular traffic in highways, see~\cite{[PH71]}. In this modeling framework, individual entities are the driver-vehicle subsystem, while the overall state of the system is described by the one particle distribution function over the microscale state of the subsystems, i.e. position and velocity. Similar to the Boltzmann equation, the dynamics of the distribution function is obtained by a balance of the number of subsystems in the elementary volume  of the space of the microscopic states, in which the inlet and outlet of particles is determined by the interactions.

Prigogine's  scientific contributions have been studied by other authors. We mention, as examples to be considered for further developments towards the modeling of behavioral systems, the introduction of the concept of heterogeneity proposed in~\cite{[PF75]} and the Enskog-type modeling of the dimension of vehicles, see~\cite{[KW1997],[KW1999A],[KW1999B]}. An important contribution to the study of multi-particle systems in traffic flow has been given by Helbing~\cite{[HEL01]}, by methods somewhat referred to statistical physics. An overview of the bibliography on vehicular traffic modeling can be found in the review~\cite{[ALBI19]}, where the recent literature on traffic, crowds and swarms is reported.

The first steps toward the development of kinetic theory methods in other fields of science were taken in the past century in the study of social dynamics in insect populations, see~\cite{[JS92]}, and in the dynamics of the competition between tumor and immune cells, see~\cite{[BF94]}. The key concepts proposed in these two papers are inspired by Prigogine's work, but also introduce new concepts. The development essentially considers the dynamics of active, rather than classical particles, where interactions are described by theoretical tools of game theory. This approach transfers the dynamics of interactions at the small scale of individuals  into the collective dynamics of populations of individuals.

This path leaves very little of the classical kinetic theory, where interactions are local, number conservative, and reversible. On the other hand, interactions involving active particles are far more complex, being non-local, non-reversible, since the output is not induced by nonlinearly additive actions. These features arise from the fact that active particles are living systems. Therefore the derivation of a mathematical theory should refer to the specificity of living systems, see~\cite{[CLE19]}. These mathematical studies interact with our society because we live in a complex world, see~\cite{[Ball12]}. The literature on physical and mathematical studies of the living matter is made of several papers and books as reported in~\cite{[VW25]}, but a unified approach is still a mythical objective.

The first attempt to organize the above pioneering literature is witnessed in the book~\cite{[BLS00]} in which various models derived within the framework of the kinetic theory, are described and critically analyzed. In particular, Chapters 5 and 6 in~\cite{[BLS00]} reports about an intense literature on the topic studied in~\cite{[BF94]}, see also~\cite{[FGP99]}; Chapter 8 is focused on the modeling of vehicular traffic, while Chapter 10, introduces the mathematical structures for more general living systems. This structure has been further modified by various authors with focus on specific applications.

The literature cited above, especially~\cite{[PH71]} and~\cite{[BF94],[BLS00]}, has motivated an intense scientific activity in this century. Various approaches have been proposed. All of them have in common that the state of the system is described by the one particle distribution function of the microscopic state of the interacting particles. This function is the dependent variable of a differential system. However, there are important differences appear in the derivation of the mathematical structures to describe the dynamics of the dependent variable, as well as in the modeling of the interactions. This research activity has led to a progressive development of the mathematical approach, but not yet to a unified mathematical theory.

The motivation for this essay, which consists mainly in the search for a mathematics of living systems through the appropriate development of methods inspired by the theoretical tools of statistical physics and, in particular, by the kinetic theory of classical particles, is presented. We then review and critically analyze the sequential steps that have developed from the first pioneering research work proposed in the last century to an intense research activity produced in this present century. Finally, we consider the search for a unified mathematical theory to describe the dynamics of living systems constituted by a large number of interacting behavioral entities. This emerging topic has been treated in books that focus on modeling, analytical topics, and  applications, see~\cite{[BEL08],[BBGO17],[BLS00],[HEL10],[PT13]}. These books present technically different approaches and  provide an important reference framework for our survey.

In particular, the book~\cite{[BBGO17]} has shown that active particle methods have been derived to describe the collective dynamics of living entities in terms of differential systems. The mathematical approach initially proposed in~\cite{[BEL08]} has been gradually developed to take into account the hints coming from applications. Although positive results have been obtained in specific case studies such as behavioral crowd dynamics, see~\cite{[BLQRS23]}, and the study of epidemics by mutating viruses, see~\cite{[B3EPT]}, we cannot naively claim that a complete theory has been obtained. Indeed, further developments are necessary, as discussed in the last section of this essay.

Therefore, the literature reviewed in this paper is mainly focused on the milestones that have signed the developments of the theory. In particular:

\vskip.1cm \noindent (i) The introduction of proliferative/destructive interactions;

\vskip.1cm  \noindent (ii) Transitions across groups of interest;

\vskip.1cm  \noindent (iii) New ideas beyond binary interactions and modeling of non-local interactions;

\vskip.1cm  \noindent (iv)  Dynamics by which individuals learn from the collectivity;

\vskip.1cm  \noindent (v) Decision dynamics by which individuals use what they have learned to apply their own strategy.

\vskip.1cm
 The outline of the contents of our paper is as follows:

\vskip.2cm \noindent Section 2 refers to a large system of interacting living entities, called  \textit{active particles}, or \textit{a-particles} for short. Their properties are defined in terms of their differences from the classical particles of the kinetic theory. These important differences  also affect the interactions that can generate proliferative-destructive events and post-Darwinian evolution and selection. Furthermore, we define the mathematical representation of the system and we define the philosophy that can lead to a mathematics for the above defined class of living systems.

\vskip.2cm  \noindent The next two sections show how the mathematical theory reviewed in the previous section can be used to derive specific models. The first section considers models with spatial dynamics, focusing on the dynamics of crowds and swarms. The second section then considers models where the dynamics in space is not taken into account. However, a richer variety of applications has been developed thanks to the modeling of interactions that include proliferative and destructive events. The presentation includes a brief review of developments in mathematical theory and advanced applications. The presentation focuses mainly on the very recent bibliography, so that we avoid repetition with respect to previous review papers.

\vskip.2cm \noindent Section 3 shows how the philosophy and methodological approach proposed in the previous section can lead first to mathematical structures thought to describe the dynamics of living systems, and then to specific models. We focus on spatial dynamics in which each a-particle has a visibility (geometrical)  and a sensitivity domain and perceives interactions within an appropriate intersection of these two domains, see~\cite{[BH17]}.  Interactions may involve a fixed number of living entities as in the conjecture proposed in~\cite{[BCC08]}. The content describes how the developments in the mathematical approach have been  driven  by  applications.

\vskip.2cm \noindent Section 4 deals with the study of models in the case of spatial homogeneity, which is considered as a special case of the dynamics treated in the previous section. The study considers models with proliferative and/or destructive encounters. Applications have been mainly focused on the immune competition, social dynamics and economics. Similar to the previous section, we show the evolution of the mathematical theory developed thanks to the motivations of applications. We have chosen immune competition as a case study to show how the theory developed.

\vskip.2cm \noindent Section 5 briefly describes the  mathematical tools of the active particle that have been developed in parallel with the methods reviewed in the previous sections. In particular,  we report on mean-field games as an important contribution to the modeling of interactions in large systems of living entities~\cite{[LL07]}; the derivation of alternative mathematical structures, such as the Fokker-Plank-Boltzmann description, see~\cite{[PT13]}; and the mathematical theory of behavioral swarms, in which  a pseudo-Newtonian dynamics is used to model the effect of interactions on the individuals of the system from~\cite{[BS12]} to~\cite{[BHLY24]}. The presentation is proposed to complete the conceptual framework to which the kinetic theory of active particles refers.

\vskip.2cm \noindent Section 6 deals with research perspectives for the further development of the kinetic theory of active particles, as well as conceivable modeling perspectives and applications. First, we consider the further development of the mathematical theory of Herbert A. Simon's philosophical theory of the artificial world, see~\cite{[Simon2019]}. We then provide an overview of the path from collective learning dynamics, see~\cite{[BDG16]}, to decision making with some considerations of multiple strategy dynamics,  see~\cite{[CG15]} for decision dynamics in complex environments, and the celebrated book by Kahnneman for the \textit{thinking} dynamics~\cite{[KA12]}. The study of these topics naturally leads to perspective ideas on artificial intelligence, see~\cite{[BDL24]}.

\vskip.2cm \noindent Section 7 presents further critical analysis and a forward look to the challenging perspective focused on the stone guest of this paper, that it is the search for a mathematics of living systems, in which the key passage is the analysis of the complexity of the collective dynamics of life. Three key topics are selected as those that can contribute to a mathematical theory of living systems. We consider the evolution of learning dynamics, the interaction of multiple strategies related to a dynamical utility function, and multiscale methods as an essential view of the dynamics of living systems. These three dynamics currently interact, while an interactive interpretation can contribute to the  mythical target mentioned above.

\section{From classical to active particles}\label{Sec:2}

We consider the collective dynamics of a large system of interacting  a-particles. This section defines their main properties and their differences with respect to the classical particles of the kinetic theory. The differences are related to the dynamics of interactions that can generate proliferative-destructive events and post-Darwinian evolution and selection. Furthermore, we define the philosophy that can lead to a mathematics for the above defined class of living systems. A critical analysis contributes to further development of the theory.

\subsection{Active particles and representations}\label{Sec:2.1}

We consider a-particles, in a plane motion, as entities whose individual state is in the following, where we use dimensionless variables.

\vskip.1cm \noindent  \textit{Position}: $\bx \in \Sigma \subseteq \RR^2$, where $\Sigma$ is the domain in which the particle can move. If $\ell$ is a characteristic length of $\Sigma$, the components of $\bx$ are divided by $\ell$. Boundaries of $\Sigma$ are denoted by  $\partial \Sigma$.

\vskip.2cm \noindent  \textit{Velocity}: $\bv = v \, \bomega$, with $\bv \in D_\bv$, where $v$ is the velocity modulus, called \textit{speed}, and $\bomega$ is a unit vector that defines the direction of the velocity. The speed is divided by $v_L$, which is the limit speed in each considered system.

\vskip.2cm \noindent  \textit{Visibility (geometrical)  and sensitivity domains} are the domains $\Omega_v$ and $\Omega_s$ in which an a-particle with velocity $\bv$ can see and feel other particles, respectively. Technical notes on the calculation of $\Omega_v$ and $\Omega_s$ are given in Section 3.

\vskip.2cm \noindent  \textit{Activity} $u \in D_u$ is the behavioral variable, called \textit{activity}, that defines the individual social-emotional state of the a-particle.

\vskip.2cm \noindent  \textit{Functional Subsystems}, for short FSs, denote the groups of interest, in which a-particles aggregate. FSs can be labeled  by $i=1, \ldots, n$.

\vskip.2cm \noindent  \textbf{Remark 2.1.} \textit{In general, the activity can be a vector. For the time being, calculations will be done with scalar activity, waiting for Section 3, where the special case of vector activity variables will be treated. For scalar variables, we will set $D_u = [0,1]$, where $u=0$ and $u=1$ correspond to the minimum and maximum values of the activity, respectively.}

\vskip.2cm \noindent  \textbf{Remark 2.2.} \textit{The spatially homogeneous case, or simply systems for which dynamics in space is not relevant, corresponds to distribution functions homogeneously distributed in space, $f_i = f_i(t, u)$. The localization can take place in a network of nodes denoted by $\kappa = 1, \ldots, K$. Then the overall state of the system is defined by $f_{i\kappa} = f_{i\kappa}(t, u)$, see~\cite{[KNO13]}.}

\vskip.2cm \noindent  \textbf{Remark 2.3.} \textit{Some specific case studies may require modeling the physical quality of the venue. Then a dimensionless parameter $\alpha \in [0,1]$ can be used where $\alpha = 0$ corresponds to very low quality and $\alpha = 1$ corresponds to very high quality. The quality of the venue increases the speed at which the a-particles move.}

\vskip.2cm

The representation of the system, according to the kinetic theory approach is delivered by the distribution function:
\begin{equation}\label{dfunction}
f_i = f_i(t, \bx, \bv, u), \hskip.5cm [0,T] \times \Sigma  \times D_\bv \times [0,1] \to \RR_+,
\end{equation}
where the term \textit{distribution function} is used to take into account that the number $n_i = n_i(t)$ of a-particles in each i-th FS depend on time due both to proliferative/desctructive interactions and dynamics across FSs.

In particular, $n_i$ is computed by zeroth order moment of $f_i$:
\begin{equation}\label{density}
n_i = n_i(t) = \int_{\Sigma \times D_\bv \times [0,1]} f_i(t, \bx, \bv, u)\, d\bx\, d\bv\, du.
\end{equation}

Note that the \textit{microscale} corresponds to the individual dynamics of a-particles. The microscopic state of a-particles is defined by position, velocity and activity, while the \textit{independent variables} are time and space. The \textit{dependent variables} are  the distribution functions $f_i$, which evolve in time and space. The \textit{macroscopic variables} are computed moments weighted by velocity and space.
\begin{equation}\label{dmoments}
M_i^{hk} = M_i^{hk} (t, \bx) =  \int_{D_\bv \times [0,1]} \bv^h \, u^k f_i(t, \bx, \bv, u)\,  d\bv\, du.
\end{equation}

\subsection{On a philosophy towards collective dynamics}\label{Sec:2.2}

Let us consider the rationale that can lead to a differential system that can describe the collective dynamics consistent with the description and representation proposed in the previous subsection. We follow the idea that the complexity of living systems and that the mathematical approach should look for differential structures that should capture the complexity features mentioned above.  Then, within this framework, methods inspired  somewhat to statistical physics can be developed to transfer the interactions at the microscopic scale into the collective dynamics described by the distribution function defined in Eq.(\ref{dfunction}).

The  concepts  outlined above will be detailed in the following, which can lead  to our \textit{philosophy}. This approach was invented and developed to tackle the difficulties of the mathematical description of the dynamics of living systems. It is worth understanding the contribution of outstanding scientists to define these difficulties.

\begin{itemize}

\vskip.2cm \item Erwin Schr\"odinger initiated the complex approach towards a mathematics for living systems, see~\cite{[ES1944]}, in which he made outstanding inventions on the way to find a correlation between physics and living systems. His contribution is not limited to an historical value. Indeed, he introduces concepts, such as  \textit{multiscale vision} and  \textit{systems biology}, which are still valid and teaches to mathematicians. Schr\"odinger developed some pioneering ideas on a systems approach in biology, where the dynamics of cells is driven by the dynamics at the molecular scale, motivated by the study of mutations (some of which are also induced by external actions such as radiation). This concept is now a fundamental point of reference the interactions between mathematics and biology, where understanding the link between the dynamics at the molecular scale of genes and the functions expressed at the level of cells is the key strategy for deriving a bio-mathematical theory.

\vskip.2cm \item Lee Hartwell has made evidence that mathematical tools used to model systems of the inert matter cannot be applied to the living matter. In fact, living entities have a purpose and chase it through behavioral strategies that are heterogeneously distributed in populations~\cite{[HART99]}. This paper has motivated the concept of behavioral variables, hereafter called \textit{activity}.
\vskip.1cm
\begin{quote}
\textit{Biological systems are very different from the physical or chemical systems of the inanimate matter.
In fact, although living systems obey the laws of physics and chemistry, the notion of function or purpose differentiate biology
 from other natural sciences. Indeed, cells are not molecules, but have a living dynamics induced by the lower scale of genes and is organized into organs.}
\end{quote}

\vskip.2cm \item Miguel A.~Herrero, Robert May,  and Robert Reed focus on the difficulty of identifying causality principles for living systems, see~\cite{[Herrero],[MAY],[Reed]}. The considerations of these scientists support the idea that causality principles, applicable to the inanimate matter, cannot be  rapidly applied to the study of living systems.

\vskip.2cm \item Ernst Mayr teaches us that mathematics should consider evolution and selection as key features of living systems, see~\cite{[MAYR]}. As a matter of fact evolution is a key feature, not only of individual entities, but also of the whole society.

\vskip.2cm \item Herbert A. Simon proposed the theory of the  \textit{Artificial World} and  \textit{Artifacts}, which teaches us that interaction dynamics depend on rules that evolve in time within an artificial environment created by all agents living in such world and modified also by external actions, see~\cite{[Simon1976],[Simon2019]}. See also the interpretation proposed in~\cite{[BE24]}.

\end{itemize}

The key feature of our approach is the assessment of the  common features of living systems to be considered  as sources of complexity. This is an open problem that deserves a deep analysis. We refer to the five features proposed in~\cite{[BBGO17]} and try to provide a deeper interpretation rather than adding new ones. In detail:

\begin{enumerate}

\vskip.2truecm  \item \textbf{Ability to express a strategy:} Living entities, i.e. \textit{a-particles}, are capable of developing specific \textit{strategies} and \textit{organizational abilities} that depend on the overall state of the surrounding a-particles and environment. These strategies evolve in time and space within the artificial world according to the philosophy of Herbert A. Simon.

\vskip.2truecm  \item \textbf{Heterogeneity:} The ability to express a strategy is not the same for all a-particles. Indeed, the \textit{expression of heterogeneous behaviors} is a common feature of a great part of living systems. Heterogeneity plays an important role in the dynamics of competition and selection, in which competition generates proliferative and destructive dynamics, while selection enforces some a-particles, and weakens others.

\vskip.2truecm   \item \textbf{Nonlinear Interactions:} Interactions are nonlinearly additive and involve immediate neighboring a-particles, but in some cases distant a-particles as well. In general, interactions are not conservative of the variables in the microscopic state of the interacting a-particles. \textit{It is important to distinguish between binary and multiple interactions. In fact, binary interactions characterize a dilute gas, but multiple interactions are very common in the case of living systems}.

\vskip.2truecm  \item \textbf{Learning Ability:} Living systems receive input from their environment and have the ability to learn from past experience by interacting with other a-particles. In most cases, learning is followed by a decision-making dynamic that determines the overall collective dynamics of the class of systems under consideration. The amount of learned information evolves over time, and consequently the interaction rules are also time-dependent.

\vskip.2truecm  \item \textbf{Darwinian Mutation and Selection:} Birth processes can produce a-particles that are better adapted to the environment, which in turn produce new entities that are better adapted to the external environment. These mutations generate selection, which also depends on the state of the a-particles.

\end{enumerate}

\vskip.2cm \noindent  \textbf{Remark 2.4.} \textit{The various mathematical approaches known in the literature consider only some of the points mentioned above. On the other hand, all of them are important, and even more could be considered. Therefore, the content of the following sections will often return to the critical considerations mentioned above.}
\vskip.2cm

Keeping all the above in mind, we can consider the approach to describe the collective dynamics of  interacting a-particles. The rationale  is shown in the flowchart in Fig.~1, where the contents in each block are as follows:

\vskip.2cm \noindent \textbf{A:} The first step is to derive of a general differential framework suitable for capturing the complexity features of living systems. This system should be made applicable to different systems of interacting living entities from cells to humans.

\vskip.2cm \noindent \textbf{S:} The focus is on the specific system under consideration. The task is to extract the key variables involved in the dynamics of such system, to understand the specific dynamics of interactions.

\vskip.2cm \noindent \textbf{B:} The task is to use mathematical models to describe how a-particles learn by interacting with other a-particles. The general should be specialized for each case study.  Inputs given in the blocks \textbf{C} and \textbf{D} help to describe the dynamics, see~\cite{[BDG16],[BDG16A]}, see also~\cite{[BD19]}.

\vskip.2cm \noindent \textbf{C} and \textbf{D:} Each a-particle is a carrier of subsystems that are capable of developing the decision process by elaborating what they have learned through different types of multiscale interactions, in which the a-particles correspond to the microscopic scale, while the above-mentioned subsystems correspond to the submicroscopic scale, i.e. a scale lower than that of individuals.

\vskip.2cm \noindent \textbf{E:} The information obtained in the previous step is the input to a differential system that models the dynamics by which individuals (in general entities) develop a decision process consistent with the case study under consideration.

\vskip.2cm \noindent \textbf{F:} The platform puts into practice the results obtained in the previous steps. It can be used to train technicians to explore the emerging dynamics of the real system for different choices of parameters. It can also be used to predict the emerging dynamics after calibrating the differential model with empirical data.

\vskip.2cm \noindent  \textbf{Remark 2.5.} \textit{Both learning dynamics and subsequent decision making should take into account the complex ``thinking'' dynamics of human brains, see~\cite{[KA12]}, and the connection between brain and mind through, see~\cite{[NBZ25]}, which cannot avoid the study of the biology of the ``brain system'', arguably the most complex system in the human body.}
\vskip.2cm

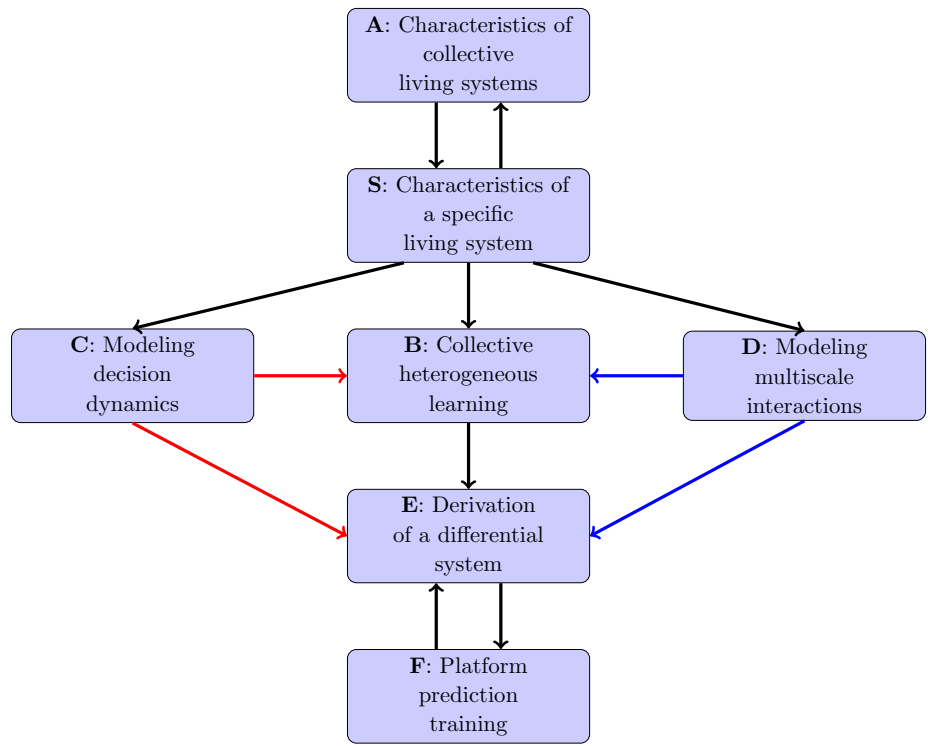
\begin{figure}[t!]
\begin{center}\scalebox{0.85}{
\begin{tikzpicture}[node distance = 2.5cm, auto]
\node [block_long] (complex) {\textbf{A}: Characteristics of\\ collective \\ living systems};
\node [block_long, below of=complex] (system) {\textbf{S}: Characteristics of\\ a specific  \\ living system};
 \node [block_long, below of=system] (learning) {\textbf{B}: Collective \\heterogeneous \\ learning};
  \node [block_long, left of=learning,xshift=-2.7cm] (micro) {\textbf{C}: Modeling \\ decision \\ dynamics};
  \node [block_long, right of=learning,xshift=2.7cm] (social) {\textbf{D}: Modeling \\multiscale\\ interactions};
   \node [block_long, below of=learning] (framework) {\textbf{E}: Derivation  \\ of  a differential \\ system};
      \node [block_long, below of=framework] (artificial){\textbf{F}: Platform \\ prediction  \\  training};
   \draw [line width=.5mm, ->] (system.south) -- (learning.north);
\draw [line width=.5mm, ->] ([xshift=-0.5cm] complex.south) -- ([xshift=-0.5cm] system.north);
 \draw [line width=.5mm, ->] ([xshift=0.5cm] system.north) -- ([xshift=0.5cm] complex.south);
  \draw [line width=.5mm, red, ->] (micro.east) -- (learning.west);
 \draw [line width=.5mm, blue, ->] (social.west) -- (learning.east);
 \draw [line width=.5mm, ->] ([xshift=-1cm] system.south) -- ( micro.north);
  \draw [line width=.5mm, ->] ([xshift=1cm] system.south) -- (social.north);
  \draw [line width=.5mm, ->] (learning.south) -- (framework.north);
 \draw [line width=.5mm, blue, ->] (social.south) -- (framework.east);
 \draw [line width=.5mm, red, ->] (micro.south) -- (framework.west);
  \draw [line width=.5mm, ->] ([xshift=-0.5cm] artificial.north) -- ([xshift=-0.5cm] framework.south);
 \draw [line width=.5mm, ->] ([xshift=0.5cm] framework.south) -- ([xshift=0.5cm] artificial.north);
\end{tikzpicture}}
\end{center}
\begin{center}
\caption{Flow chart for the overall dynamics}\label{fig1}
\end{center}
\end{figure}

\subsection{Critical analysis}

The philosophy presented in the previous subsection leads to a strategy suitable for deriving mathematical models to describe the collective dynamics of living, i.e. behavioral, systems. Applications have been developed in various fields of science. The following sections report on a bibliography selected according to the framework defined in Section 1. These applications have shown that mathematical models can describe important features of collective dynamics. These encouraging results have motivated further studies and have raised a key question that deserves an answer.
\vskip.1cm
\begin{quote}
\textbf{Key question:} \textit{To what extent can the mathematical tools to be derived according to the above general framework be considered a mathematical theory?}
\end{quote}
\vskip.1cm

The answer can be given in the last section, after the derivation of the mathematical approach and an overview and critical analysis of the literature in this field.  As a first step, we specify that we study the collective dynamics of interacting living entities. This implies that we have to consider complex systems and complexity. David C.~Krakauer's book provides an excellent historical and philosophical reference, see~\cite{[Krakauer]}. Furthermore, a definition of \textit{complexity} is useful, as already done in~\cite{[BBD21]}, we extract from~\cite{[Simon1962]} (in the first lines of page 468) the following statement:
\vskip.1cm
\begin{quote}
 \textit{Roughly speaking, by a complex system I mean one that consists of a large number of parts that interact in non-simple ways. In such systems, the whole is more than the sum of the parts, not in some ultimate metaphysical sense, but in the important pragmatic sense that it is not trivial to infer the properties of the whole from the properties of the parts and the laws of their interaction.}
\end{quote}
\vskip.1cm

In fact, the paper cited above teaches us much more than a definition of complexity. In fact, it points out the conceptual differences between systems with spatial dynamics and those limited to temporal dynamics. Moreover, it highlights the difficulty of extrapolating the dynamics of a few parts to the dynamics of the whole. In any case, the above sentence has encouraged our work.

Then, some preliminary speculations are proposed in view of the mathematical tools developed in the next sections and to the answer to the key question that we have posed above.

\begin{itemize}

\vskip.2cm \item  \textit{Modeling interactions:} Binary interactions correspond to a specific property of rarefied gases, in particular the \textit{Boltzmann equation, in which the spatial distance between particles makes the probability of multiple interactions negligible}. Indeed, the structures derived by this assumption preserve some of the properties of classical kinetic theory, but the interactions are not reversible and do not preserve mechanical quantities, see~\cite{[Kogan]} and \cite{[CIP93]}. Therefore, these considerations cannot be applied to systems of living entities. An analogous consideration can be focused on whole-field interactions corresponding to the encounters of each particle with all particles, generally with linearly additive dynamics. In this case, the reference model of the kinetic theory of classical particles is the Vlasov equation, see~\cite{[CIP93]}.

\vskip.2cm  \item  \textit{Multiple and multiscale interactions} are  somewhat inspired by the conjecture of Giorgio Parisi's team, see~\cite{[BCC08]}, which is based on the idea that each a-particle \textit{learns} from a fixed number of other a-particles within its sensitivity domain.  \textit{Interactions involve multiple entities and can be multiscale}, for example, when an a-particle also interacts with all particles as a whole represented by low-order moments of the distribution function. On the other hand, interactions can be driven by inputs from  lower scales such as from the molecular (gene) scale in cells or neuronal networks in humans.

\vskip.2cm  \item  \textit{Nonlinearity:} In general, the mathematical structures, derived in the next sections, are nonlinear in the dependent variables defined by Eq.~(\ref{dfunction}). In fact, in most cases the product of the distribution functions, involved in the interactions,  appears. On the other hand, the concept of nonlinear interaction refers to the mapping of the microscale state from the pre-interaction to  the post-interaction state. \textit{Interactions are defined as nonlinear if the map depends, in some way, on the  $f_i$}. On the other hand, if the dependence is just on the microstates of the a-particles, then \textit{the interactions can be considered linear}.

\end{itemize}

The above considerations can contribute to the study of modeling, analytical problems, and computation that will be developed in the next sections. Furthermore, they can help to understand that we are moving far away from the mathematical tools of statistical physics to new methods that are candidates for a theory. The effort to look for a new physics and, as a consequence, a new mathematics has been  the key idea of Schr\"odinger's seminal work~\cite{[ES1944]}. It was a target specifically focused on economics~\cite{[JL21]}. We have already emphasized most of the conceptual differences and, in particular, the fact that the dynamics of living systems appear far from equilibrium. A sharp interpretation, worthy of  consideration, has been proposed by A.Aristov in~\cite{[ARI19]}.

\section{Mathematical theory of active particles with spatial dynamics}\label{Sec:3}

This section shows how the general philosophy proposed in Section 2 can be translated into a mathematical tool for describing behavioral systems in general. We follow the guidelines of the mathematical theory of active particles proposed in~\cite{[BBGO17]} and further developed in~\cite{[BBD21]}.
Some further developments are proposed in this present paper. The presentation is done in three steps. First, we study the dynamics of the interactions. Then, we derive the mathematical structure we are looking for. Finally, we develop a critical analysis of the literature in the field and we show how the mathematical theory developed as motivated by the applications.

\subsection{On the concept of linear and nonlinear interactions}\label{Sec:3.1}

The above mentioned structure transfers the dynamics at the low individual scale to the collective motion described by the kinetic theory of active particles. Therefore, a key step to be taken is the study of these interactions. For simplicity of notation, we consider a system consisting of only one FS in the absence of proliferative and destructive events. The generalization to different FSs follows. To specify the interaction rules, we need to introduce three types of a-particles, which may correspond to different physical systems, such as the vehicle-driver subsystem, pedestrians, and others. In addition, we also define the interaction domain.

\vskip.2cm \noindent \textbf{Test particles} with distribution function $f(t,\bx, \bv, \bu)$, which are representative of the whole system.

\vskip.2cm \noindent  \textbf{Field particles} with distribution function $f(t, \bx^*,  \bv^*, \bu^*)$, which, by interacting with the
test particles, can lose their micro-state;

\vskip.2cm \noindent \textbf{Candidate particles} with distribution function $f(t, \bx_*,  \bv_*, \bu_*)$, which can acquire, in probability, the micro-state of the test particle after interaction with the field particles.

\vskip.2cm \noindent \textbf{Interaction domain} in which the interactions of \textit{test} and \textit{candidate} with \textit{field} particles take place. In general, it is a conical region $\Omega = \Omega (t, \bx, \bv;\theta, R)$, whose vertex is in $\bx$, whose axis being the direction of velocity $\bv$, while $\theta$ and $R$ are variables that define the width of the cone and the visibility/sensitivity distance. In this domain, test or candidate a-particles learn the state of the field particles.
\vskip.2cm
This definition is not yet precise and can be applied either to the so-called \textit{metrical interactions}, which depend on the microscopic variable but not on the distribution function, or to the \textit{topological interactions}, which depend on the local density. Detailed calculations have been given in~\cite{[BS12]}. In the following, we will generically use the notation $\Omega$ to denote the interaction domain in the derivation of the differential system, while the topic will be discussed and detailed calculations will be given later.

Interactions induce a modification of the microscopic state, i.e. activity, velocity direction, and speed. In general, according to the theory developed in~\cite{[BBGO17]}, the dynamic depends on the  micro-state and distribution function of the pedestrians in the interaction domain. The theory proposes a dynamics of interactions modeled as follows:

\begin{itemize}
\vskip.2cm \item   Interactions (within $\Omega$) occur with the \textit{interaction rate}
$$
\eta[f](\bx, \bx^*, \bv_*,\bv^*, \bu_*, \bu^*),
$$
which models the interactions between test, or candidate particles, with field particles.

\vskip.2cm \item The state of  \textit{test} and \textit{candidate} is modified according to the \textit{transition probability density}
$$
\cA[f](\bv_* \to \bv, \bu_* \to \bu|\bx, \bx^*, \bv_*,  \bv^*, \bu_*, \bu^*),
$$
which  models the probability that a candidate particle at $\bx$ with state $\{\bv_*, \bu_*\}$ shifts to the state of the test particle $\{\bv, \bu\}$ due to the interaction with field particles with state  $\{\bv^*, \bu^*\}$ in $\Omega$.

\end{itemize}

Let us now put in evidence some specific features of the interactions involving active particles, which are different from those of classical particles as clearly stated in the books~\cite{[BBGO17],[PT13]}. The following remarks focus on them, while other details are given when specific applications are critically analyzed.

\vskip.2cm \noindent  \textbf{Remark 3.1.}  \textit{Unlike those of classical particle theory, the interactions are neither reversible nor local. In fact, the active particles are sensitive to the other particles at a distance. The interactions can depend not only on the state of the interacting pairs, but also on $f$ and the specific properties of the domain in which the dynamics takes place. This dependence leads to nonlinearity of the dynamics of the interactions, as discussed in Section 2.}

\vskip.2cm \noindent  \textbf{Remark 3.2.}   \textit{The dynamics of interactions should be developed in two successive steps. The \textit{learning dynamics}is followed by the decision dynamics. Learning, by candidate or test particles, refers to both the microscopic and macroscopic state of the field particles in the sensitivity domain $\Omega$.}
\vskip.2cm

The concepts expressed in the above remarks need some additional reasoning to make precise the nonlinearity features that should be referred to the structures that provide a differential system considered to describe the collective dynamics.

\subsection{On the derivation of mathematical structures}\label{Sec:3.2}

Consider first a system with a single FS. The general mathematical structure for time evolution of the distribution function $f$ can be obtained from a balance of particles in the elementary volume of the space of the micro-states $[\bx, \bx + d\bx] \times [\bv, \bv + d\bv] \times [u, \bu + d\bu]$. This  structure is derived by equating  the rate-of-change of the number of active particles (a-particles for short) plus the transport due to the velocity variable  to the net  flux rate within the elementary volume, that is
\begin{equation}\label{eq.meso}
	\frac{\partial f}{\partial t}+\bv \cdot \nabla_\bx \, f = J[f](t, \bx, \bv, \bu)  = (\cG[f] - \cL[f])(t, \bx, \bv, \bu),
\end{equation}
where  the dot product denotes the standard inner product in $\RR^2$,  $\nabla_\bx$ denotes the gradient operator with respect to the space variables only,  $ J[f](t, \bx, \bv, \bu)$ is the net  flux rate due to interactions, which is given by the difference between the \textit{gain} and \textit{loss} terms, $\cG$ and $\cL$, both nonlinearly  acting on $f$, of pedestrians in the elementary volume of the phase space about the test microscopic state $(\bx,\,\bv, \, \bu)$. Note that, here and  in below, the square brackets denotes functional dependence with respect to their arguments, e.g., dependance on the spatial derivatives of the arguments in brackets.

Based on these concepts,  the gain and loss terms can be written as follows:
\begin{eqnarray}   \label{Gain}
 && \cG[f](t, \bx, \bv, \bu)  = \int_{\Gamma\times {D_\bv} \times {D_\bu}} \eta[f](\bx, \bx^*, \bv_*,\bv^*, \bu_*, \bu^*) \nonumber \\[3mm]
 && \hskip1cm  \times \cA[f](\bv_* \to \bv, \bu_* \to \bu|\bx, \bx^*, \bv_*,  \bv^*, \bu_*, \bu^*) \nonumber \\[3mm]
 && \hskip1cm  \times  \, f(t, \bx, \bv_*, \bu_*) f(t, \bx^*,  \bv^*, \bu^*)\, d\bx^*\,d\bv_*\,d\bv^* \, d\bu_*\, d\bu^*,
\end{eqnarray}
and
\begin{eqnarray} \label{Loss}
&& \cL[f](t, \bx, \bv, \bu) = f(t, \bx, \bv,\bu) \int_{\Gamma} \eta[f](\bx,\bx^*, \bv,\bv^*,\bu, \bu^*)\nonumber \\[3mm]
 && \hskip1truecm \times \,f(t, \bx^*, \bv^*, \bu^*)\,d\bx^* \, d\bv^* \, d\bu^*,
 \end{eqnarray}
where $\Gamma = \Omega \times {D_\bv} \times {D_\bu}$.

The macroscopic observable quantities can be defined, under suitable integrability assumptions, by weighted  moments of the distribution function. For instance, the local \textit{density}  reads
\begin{equation}\label{mac-1}
\rho (t, \bx) = \int_{D_\bv \times D_\bu} f(t,\,\bx,\,\bv,\, \bu)\,d\bv\,d\bu,
\end{equation}
where $D_\bv = [0, 2 \pi) \times [0,1]$ in polar coordinates, and the \textit{mean velocity} is defined as
\begin{equation}\label{mac-2}
\boldsymbol{\xi} (t, \bx) = \frac{1}{\rho (t,\bx)} \int_{D_\bv \times D_\bu}  \bv\, f(t,\,\bx,\,\bv,\, \bu)\,d\bv\,d\bu.
\end{equation}

These mathematical structures provide the mathematical framework to be used to derive specific models, which are obtained by using a phenomenological interpretation of interactions to derive the above mentioned terms $\eta$ and $\cA$. These differ for each case study according to the specific dynamics of interactions.

The derivation of the mathematical structure for systems made of a finite number of interacting FSs is immediate. We denote the FSs by $j=1, \ldots, n$. Then, we have to define the interaction terms as follows:

\vskip.2cm \noindent $\eta_{ij}$ is the interaction rate between a-particles of the i-th and the j-th FS.

\vskip.2cm \noindent $\cA_{ij}$ is the transition probability density corresponding to the interaction between a-particles of the i-th and the j-th FS.

\vskip.2cm  Both $\eta_{ij}$ and $\cA_{ij}$ can depend on the microscopic state and on their distribution function. Let $\bbf = \{f_1, \ldots, f_n\}$. Then, technical calculations analogous to those we have seen above lead to the following:
\begin{equation}\label{MultiFS}
\frac{\partial f_i}{\partial t}+\bv \cdot \nabla_\bx \, f_i = J_i[\bbf](t, \bx, \bv, \bu) = (\cG_i[\bbf] - \cL_i[\bbf])(t, \bx, \bv, \bu),
\end{equation}
where
\begin{eqnarray} \label{G-MultiFS}
&&\cG_i[f_i, f_j] = \sum_{j=1}^n \int_{\Gamma \times {D_\bv} \times {D_\bu}} \eta_{ij}[f_i, f_j](\bx, \bx^*, \bv_*,\bv^*, \bu_*, \bu^*)\nonumber\\[3mm]
 &&  \hskip1.5cm \times \cA_{ij}[f_i, f_j](\bv_* \to \bv, \bu_* \to \bu|\bx, \bx^*, \bv_*,  \bv^*, \bu_*, \bu^*) \nonumber \\[3mm]
&& \hskip1.5cm  \times  \, f_i(t, \bx, \bv_*, \bu_*) f_j(t, \bx^*,  \bv^*, \bu^*)\, d\bx^*\,d\bv_*\,d\bv^* \, d\bu_*\, d\bu^*,
\end{eqnarray}
and
\begin{eqnarray} \label{L-MultiFS}
&&  \cL_i[f_i, f_j] = f_i(t, \bx, \bv,\bu) \sum_{j=1}^n \int_{\Gamma} \eta_{ij}[f_i, f_j](\bx,\bx^*, \bv,\bv^*,\bu, \bu^*)\nonumber \\[3mm]
 && \hskip2truecm \times \,f_j(t, \bx^*, \bv^*, \bu^*)\,d\bx^* \, d\bv^* \, d\bu^*,
\end{eqnarray}
were $\Gamma = \Omega \times {D_\bv} \times {D_\bu}$.

\vskip.2cm \noindent  \textbf{Remark 3.3.} \textit{The structures (\ref{eq.meso}) and (\ref{MultiFS}) can include in the transport terms external actions that act on both the mechanical and the behavioral variables. This topic is covered sections (3.3.4)-(3.3.5) in the book~\cite{[BBGO17]}. The same subsection shows how a decay term can be added to the right-hand side of the equation when the behavioral variable has a natural tendency toward equilibrium. This concept has been used in models of immune competition, see~\cite{[BDK13]} for competition between cancer and immune cells, see~\cite{[BEF24]} and~\cite{[BK24]} for competition with viruses. In particular, activated immunity drops to the level of innate immunity when pathology is suppressed.}
\vskip.2cm

Some technical details clarify the statements of this remark. We refer to Eq.(\ref{eq.meso}) and consider how the transport term of this equation can be modified by two different types of action: an external field $\bF = \bF[f](t,\bx)$, per unit mass, which affects the mechanical dynamics of the a-particles but does not modify the activity, and an external action $\vf = \vf[f](t, u)$, which acts uniformly on the a-particle in space but induces a transport term for the activity variable. For simplicity, here and below the activity variable is supposed as a scalar.

In the first case, the structure is modified by adding the action of the external field $\bF[f](t,\bx)$ to the left term, as follows:
\begin{equation}\label{ext.meso}
	\frac{\partial f}{\partial t}+\bv \cdot \nabla_\bx f + \nabla_\bv \big(\bF[f] \, f\big) = J[f](t, \bx, \bv, u),
\end{equation}
where the operator $J[f]$ corresponds to the difference between the gain and loss terms defined in Eqs.{\ref{Gain}, \ref{Loss}).
In the second case, the formal structure should also include the term $\partial \big( u \vf[f] \, f\big)$ on the left side. However, the physics of the small-scale dynamics must be related to the specific case study under consideration.}

The above presentation can be completed by some considerations on the calculation of the $\Omega$ according to Parisi's conjecture, see~\cite{[BCC08]}. Some hints are proposed in the following for a system with only one functional subsystem, i.e. referring to Eq.(\ref{eq.meso}). These considerations refer to~\cite{[BS12]} and to~\cite{[BH17]} for analytic topics, see also~\cite{[ALBI24],[ALBI-P12]}.

First we observe that if a critical density $\rho_c$ can be defined (in the swarm theory of birds, it can be a fixed number of them), $\rho_c$ can be put on the left-side of Eq.~(\ref{mac-1}). Then, using polar coordinates in the integral on the right side, and supposing symmetrical interactions, i.e. an arc of circle in the plane or a conical domain in space, in both cases with  symmetrical semi-amplitude $\theta$, it is possible computing the radius $R$ of the so-called \textit{sensory domain} $\Omega_s[f]$, see Fig.~2. This is a theoretical calculation that has to be related to the \textit{visibility domain} $\Omega_v$.

Then, the \textit{effective interaction domain} $\Omega = \Omega[f]$ is defined as follows:
\begin{equation}\label{Omega}
\Omega_s \subseteq \Omega_v \hskip.2cm \Rightarrow \hskip.2cm \Omega [f] = \Omega_s, \hskip.5cm \Omega_v \subset \Omega_s \hskip.2cm \Rightarrow \hskip.2cm \Omega [f] = \Omega_v,
\end{equation}
by a combination or switching between metric  and topological domains are also possible, still with $\Omega$ depending on $f$.

The figure below shows  $\Omega_s$ and  $\Omega_v$  in a plane representation. A similar one can be done for a three dimensional representation by showing a conical volume in space.

\begin{figure}[htpb] \label{domains}
\begin{center}
\begin{tikzpicture}[scale=0.85]

\draw[rotate=30, xscale=1.9,yscale=1.3] (0.8,0) circle (1.8); 
\draw [fill=blue!10] (0.2,0)--(1.8,0.3) arc [radius=2, start angle=20, end angle= 56] -- cycle; 
\draw [violet] (0.2,0)--(2.3,0.2) arc [radius=2.4, start angle=10, end angle= 60] -- cycle;   

 \draw[fill] (0.2,0) circle [radius=2pt];
 \draw (0.2,0.2) node[anchor=east] {Particle};
  \draw (0.2,0) node[anchor=north] {$\bx$};
 \draw[blue,->] (0.2,0)-- (1.1,0.6) ;
 \draw [blue](1.1,0.6) node[anchor=south] {$\bv$};  

\draw[dashed,->] (1,-0.21) -- (1,0.4);
\draw (1,-0.21) node[anchor=north] {$\Omega_s$};

\draw[dashed,->] (1.9,-0.21) -- (1.9,0.5);
\draw (1.9,-0.21) node[anchor=north] {$\Omega_v$};

\end{tikzpicture}
\begin{tikzpicture}[scale=0.85]

\draw[rotate=30, xscale=1.9,yscale=1.3] (0.8,0) circle (1.8); 
\draw [fill=blue!10] (0.2,0)--(1.8,0.3) arc [radius=2, start angle=20, end angle= 56] -- cycle; 
\draw [violet] (0.2,0)--(2.3,0.2) arc [radius=2.4, start angle=10, end angle= 60] -- cycle;   

 \draw[fill] (0.2,0) circle [radius=2pt];
 \draw (0.2,0.2) node[anchor=east] {Particle};
  \draw (0.2,0) node[anchor=north] {$\bx$};
 \draw[blue,->] (0.2,0)-- (1.1,0.6) ;
 \draw [blue](1.1,0.6) node[anchor=south] {$\bv$};  

\draw[dashed,->] (1,-0.21) -- (1,0.4);
\draw (1,-0.21) node[anchor=north] {$\Omega_v$};

\draw[dashed,->] (1.9,-0.21) -- (1.9,0.5);
\draw (1.9,-0.21) node[anchor=north] {$\Omega_s$};

\end{tikzpicture}
\caption{Sensory, visibility, and interaction domains $\Omega_s$, $\Omega_v$, and $\Omega=\Omega[f] = \Omega_s \cap \Omega_v$.
Left: $\Omega_s \subset \Omega_v \Rightarrow \Omega = \Omega_s$. \hskip.2cm Right: $\Omega_v \subset \Omega_s \Rightarrow \Omega = \Omega_v$.}
\end{center}
\end{figure}
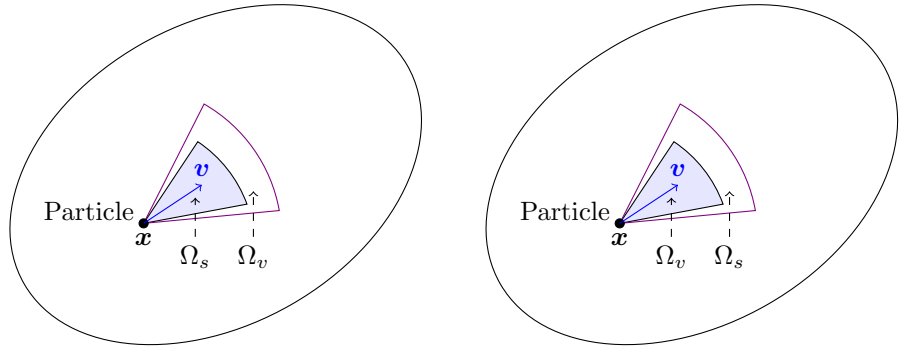

\vskip.2cm \noindent  \textbf{Remark 3.4.} \textit{The modeling of the sensitivity domain has been developed for a single functional subsystem. Therefore the notation $\Omega = \Omega [f]$ refers to Eq.(\ref{eq.meso}) and it is consistent with the above mentioned study. On the other hand, the case of interaction functional subsystems is object of ongoing studies, in which it is  to understand how different types of active particles modify $\Omega$. Known applications have used the visibility domain as sensory domain, i.e. $\Omega = \Omega_s$.}
\vskip.2cm

The sensory domain may evolve over time, each individual may have her/his own way of perceiving, or it may be a multisensory integration of auditory, olfactory, visual, or taste.  Each individual may perceive and/or integrate information from others in this sensory domain or from the  environment, for instance, the moving waves, short or long range signals of various types coming from other individuals, the temperature, electromagnetic, gravity or pressure gradients from the environment \cite{[CL06],[PS09],[PA98]}.

Then, the interaction domain can be viewed as related to the sensory domain, and as well to $f$. The interaction domain can be \textit{metric} or \textit{topological}, depending on different case studies, but no single criterion applies to all. In a {\it metric-distance domain}, an individual interacts with all other individuals in its neighborhood. The number of members in the metric domain fluctuates as nearby individuals enter and leave its neighborhood. In this model, a fixed perceptual boundary is assumed   corresponding to their interaction neighborhoods is assumed and should be specified by an interaction range and/or radius. In {\it topological-distance domains}, the individual's  neighborhood  is elastic continuously expanding and shrinking. Each individual chooses to interact with at most a given number of neighbors, regardless of the absolute distance to the reference individual. For example, the number of members in a topological domain of birds ranges between 6-7 neighbors \cite{[BCC08]}.

The interaction domain can be further subdivided into different parts, so that the interaction rules of each part are different. For example, in Reynolds' seminal paper \cite{[R87]} on flocking, the individuals of the simulated flock, called boids, adjust their dynamic behavior by following three simple rules: {\it separation} (collision avoidance from neighbors too close together), {\it cohesion} (attraction to distant neighbors), and {\it alignment} (in particular, velocity alignment with neighbors). For each boid, Reynolds calculates the average of the summed direction vectors of the boids that are in view and in sight. The final direction of a boid's movement is calculated as the weighted sum of the results of these three rules. Varying the sizes of these zones could reproduce different aggregation patterns, see for example~\cite{[CK02],[HH08]}.

The sequential steps  of the dynamics of interactions can be summarized in the  following way:

\vskip.2cm \noindent \textbf{Step 1. Sensing:} Each candidate a-particle interacts with the field a-particles in $\Omega[f_i]$. This action occurs with rate $\eta_{ij}$ and it is  specific to the microscopic state of the interacting particles.

\vskip.2cm \noindent  \textbf{Step 2. Learning:} Each candidate a-particle learns the microscopic states  of the field a-particles in $\Omega[f_i]$. This is a collective learning dynamic that has to be properly modeled for each case study under consideration.

\vskip.2cm \noindent  \textbf{Step 3. Decision-Making:} For the gain dynamics, each candidate a-particle, after the learning mentioned in the previous item, changes the microscopic state into a new one with probability by a transition modeled by  $\cA_{ij}$. For the loss dynamics, test a-particles undergo analogous interactions with field a-particles in $\Omega[f_i]$ and lose their microscopic state.

\vskip.2cm \noindent  \textbf{Step 4. Collective dynamics:} It is the consequence of the previous sequential steps. Such dynamic is obtained by the differential systems reported in Eq.~(\ref{eq.meso}), with gain and loss terms defined in  Eq.~(\ref{Gain}) and  Eq.~(\ref{Loss}) or by the analogous ones for multi-component flows.

\vskip.2cm

The above sequential steps have been applied in various applications, in particular traffic, crowds and swarms. The following remarks introduce to the applications selected and reviewed in the next section. Applications need to understand how the sensory domain $\Omega$ can be related to the visibility domain $\Omega_v$. Some preliminary considerations are reported in the following remark.

\vskip.2cm \noindent  \textbf{Remark 3.5.} \textit{Candidate a-particles receive a full information if $\Omega \subseteq \Omega_v$. In this case, the decision making described by $\cA[f]$ (or $\cA_{ij}[f_i, f_j]$), can be obtained by  interpreting the way of thinking of the interacting individuals. On the other hand,  \textit{if $\Omega_v \subset \Omega$, the decision-making is modified  with respect to the above mentioned way of thinking}. A mathematical description of this modification is a challenging objective. A simple assumption (conjecture) is that $\cA[f]$ can only be used for integration over $\Omega_v$,  but the output is random for the integration over the remaining part of $\Omega$. }

\subsection{An overview  of the literature and critical analysis}\label{Sec:3.3}

The structure in Eq.~(\ref{eq.meso}), with the gain and loss terms defined in Eq.~(\ref{Gain}) and Eq.~(\ref{Loss}), corresponds to the current state of the art. It is worth understanding how the mathematical theory developed from the seminal papers~\cite{[PH71]} and~\cite{[BF94]}. The following brief review focuses on how the mathematical approach reached the current state of the art. The presentation does not claim to be complete, it has been selected according to the research experience of the authors, it may also be subject to individual bias. Specifically, our selection focuses on research that is developed in the search of a mathematics of living systems, rather than technical applications of known techniques.

As mentioned above, the first indication of the use of methods of kinetic theory is the modeling of vehicle dynamics in highways. Prigogine's work was enriched by the concepts of heterogeneity proposed in~\cite{[PF75]}, while the paper~\cite{[BF94]} on immune competition showed how the microscopic state of interacting entities can be described by a behavioral variable. We should also mention the work of the Kaiserslautern team focused on real-world applications, see~\cite{[KW1997],[KW1999A],[KW1999B]}. The developments were mainly motivated by applications, mainly devoted to the study of vehicular traffic, crowds, swarms, and multicellular systems.

A study of vehicular traffic introduced the idea that interactions depend not only on the state of the interacting pairs, but also on the local density, see~\cite{[CDF07],[DT07]}. This feature has been exploited in the study of interactions between aggregations of fast and slow vehicles. The onset of fragmented flow is shown by simulation and confirmed by empirical data, see~\cite{[BDF12]}. Models with a finite number of speeds have been introduced in~\cite{[FT13]}, while multilane dynamics has been studied in~\cite{[Zagour23]}.

Applications of the kinetic theory of active particles to the study of crowd dynamics have already been reviewed in~\cite{[BLQRS23]}, see also the review~\cite{[HAG21]}. Therefore, we will limit ourselves to describing how the approach to this topic has evolved in the last decade. In fact, these considerations may also help to understand how the theory has developed over time.

\begin{enumerate}
\item \vskip.2cm A mathematical model for a differential system with a discrete number of velocities has been proposed in~\cite{[BBK13]}. This technicality leads to a realistic model of the choice of velocity direction, and it also allows the dynamics to be described by a system of PDEs rather than by the integro-differential system, see~\cite{[KQ19]}. The dynamics of attraction-repulsion has been studied in~\cite{[TM25]}, while the role of  asymmetric interactions on the collective dynamics has been developed in~\cite{[GMS19]}. The study of human crowds has often been related to safety problems, see~\cite{[HS17],[RNC16],[RN18],[WCV16]}, also referred to as the psychology of crowds, see~\cite{[WSB17]}, which may be influenced by specific local environments. The role of leaders can be modeled to improve safety conditions, see~\cite{[AFS22]}.

\vskip.2cm \item  Modeling the emotional state by using a specific parameter to model the level of stress has been proposed in~\cite{[BG15]} and applied to study the dynamics of interacting, even antagonistic, groups, see~\cite{[BGO19]}. This approach has been further developed to account for different types of behavioral states, such as contagion awareness, see~\cite{[KQ20]} and~\cite{[ABKT23],[BEH25],[KOOQ21]}. Often, the focus is related to congestion and safety, see~\cite{[LLH20],[ZXW16]}, while the study of interactions between crowds and mobile structures has been developed in~\cite{[VB09]}.

\vskip.2cm \item The  decision takes into account the following tendencies: search for the programmed direction, search for less crowded areas and attraction from the main stream. As observed in~\cite{[BDL24]}, the decision is obtained by a weighted sum of the three tendencies, where the weights depend on the local density and on the level of stress. In particular, increasing density increases the search for the vacuum direction, while increasing activity (in the case of stress) increases the attraction to the  main stream,  while once  the direction is chosen, pedestrians adjust their speed to the local density conditions.

\end{enumerate}

The contents of this section have provided an overview and critical analysis of the modeling of a-particles using kinetic theory methods.  However, additional topics can be briefly reported looking ahead to research perspectives. The following notes have also been selected according to the research experience of the authors.

\vskip.2cm \noindent $\bullet$ \textbf{Swarms and self-propelled particles:} Considerations analogous to those reported for crowd dynamics can be addressed to the kinetic theory approach to swarm dynamics, see~\cite{[BS12]} and~\cite{[BH17]}. The main difference, with respect to swarm dynamics, is in the modeling of interactions, in fact, trend attraction is predominant with respect to other possible trends. The classical swarm theory, formulated in the famous paper by Cucker and Smale~\cite{[CS07]}, is based on the assumption of a homogeneous distribution of the activity variable. However, the development of a new theory, where the activity variable is an internal variable distributed over the whole population, was promoted in~\cite{[BHO20]} and further developed in some subsequent papers, which will be properly discussed in Section 5, focusing on mathematical alternatives to the kinetic theory of active particles.

\vskip.2cm \noindent $\bullet$ \textbf{Dynamics of multicellular systems:} Applications to cellular dynamics could be developed with reference to~\cite{[CDS23],[CS21],[CKS21]}. The mathematical structure is analogous to that presented in this section, but the decision making takes into account the trends generated by the biological functions that define the \textit{activity} for each biological dynamics object of the modeling approach. An additional feature to take into account is that the motion of cells evolves through interactions with the extracellular matrix generated by the cells themselves. Therefore, the modeling approach should take into account the biological functions, i.e. activities, expressed by the cell populations, in particular the FSs.  Biological activities can induce proliferative and destructive dynamics. Another interesting application is the modeling of multiple sclerosis using kinetic and reaction diffusion equations, see~\cite{[BGMT24],[OT24]}.

\vskip.2cm \noindent $\bullet$ \textbf{Exogenous networks:} The dynamics in space can be modeled by exogenous networks connecting nodes localized in space. In this approach, the dynamics are modeled by migrations across nodes, see~\cite{[KNO13]}. These methods have been applied to the  spread of opinions, see~\cite{[KNO13]}, and epidemics, see~\cite{[ADKV]}. Arguably, further developments can go beyond population dynamics and focus on the modeling of spatial dynamics in biological environments, such as in-host dynamics, see~\cite{[BEF24]}. New ideas on the dynamics over networks have been developed in~\cite{[BDZ25]}, takes into account the heterogeneous properties of the network, which is a very important feature to be considered in heterogeneous systems.

\vskip.2cm \noindent $\bullet$ \textbf{Multiscale dynamics:} The above content has been focused mainly on mathematical tools at the mesoscopic (kinetic) scale. However, a multiscale interpretation is important, even necessary. In fact, the dynamics at the microscopic  (individual-based) scale, which is generally modeled by ordinary differential equations or agent methods contribute to modeling interactions in the kinetic equations.
\vskip.2cm
 The study of multiscale dynamics has posed new problems to mathematics. Let us emphasize the following  points, see Figure 3.

\vskip.2cm \noindent  \textbf{Remark 3.6.} \textit{Microscopic models can be derived by studying interactions at the microscopic scale, taking into account nonlinearities and nonlocal dynamics. Microscopic models contribute to the derivation of kinetic-type equations. The decision making by which individuals choose their trajectories is generated by brain reasoning developed at the sub-microscopic scale of neurons. Finally, the derivation of macroscopic models can be obtained by asymptotic methods according to the clues of the Sixth Hilbert Problem, see~\cite{[Hilbert]}.
This vision requires that the derivation at each scale is based on the same principles, this topic has been developed with a focus on crowd dynamics, see~\cite{[BGQR22]}.}
\vskip.2cm

The literature on the micro-macro derivation presents highly challenging difficulties related to the mathematical structure of the Boltzmann equation, which presents conceptual difficulties, as shown in in~\cite{[DHX25]}, see also~\cite{[DH23],[DH25]}. The conceptual difficulty is somewhat related to loss or regularity for long time behaviors shown in the celebrated paper by DiPerna and Lions, see~\cite{[DL1989]}.

The derivation of macroscopic models for active particle methods can be developed by asymptotic, Hilbert, methods, with a focus on the propagation of epidemics, see~\cite{[BC22],[BC23],[ZS22]}. See also~\cite{[BKMTZ24],[CKS21],[ZS22]}.  An important review has focused on the derivation of cross-diffusion models from the underlying description at the scale of cells, see~\cite{[BKZ24]}. The authors support the analytical approach with simulations showing the dynamics of pattern formation.

 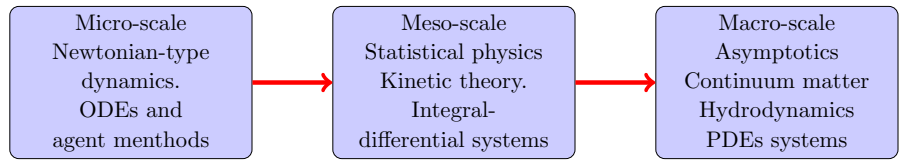
\begin{figure}
\tikzstyle{decision} = [diamond, draw, fill=blue!20,
    text width=6.0em, text badly centered, node distance=5cm, inner sep=0pt]
\tikzstyle{block1} = [rectangle, draw, fill=green!20,
    text width=6.0em, text centered, rounded corners, minimum height=4em]
    \tikzstyle{block3} = [rectangle, draw, fill=red!20,
    text width=6.0em, text centered, rounded corners, minimum height=4em]
    \tikzstyle{block2} = [rectangle, draw, fill=blue!60,
    text width=6.0em, text centered, rounded corners, minimum height=4em]
    \tikzstyle{block7} = [rectangle, draw, fill=white!60,
    text width=6.0em, text centered, rounded corners, minimum height=4em]
\tikzstyle{line} = [draw, -latex']
\tikzstyle{cloud} = [draw, ellipse,fill=red!20, node distance=3cm,
    minimum height=2em]
    \begin{center}\scalebox{.85}{
\begin{tikzpicture}[node distance = 5cm, auto]
    \node [block_long] (Inf*j) {Meso-scale\\Statistical physics\\ Kinetic theory.\\Integral-differential systems};
    \node [block_long, left of=Inf*j] (micro) {Micro-scale\\ Newtonian-type\\  dynamics.\\ ODEs and agent menthods};
    \node [block_long, right of=Inf*j] (Inf*j+1) {Macro-scale\\ Asymptotics\\ Continuum matter\\ Hydrodynamics \\ PDEs systems};
    \draw [line width=.7mm, red, -> ] (micro) -- node  {} (Inf*j);
    \draw [line width=.7mm, red, -> ] (Inf*j) -- node  {} (Inf*j+1);
    \end{tikzpicture}}
    \end{center}
    \begin{center}
\caption{Multiscale dynamics}\label{fig3}
\end{center}
\end{figure}

\section{Mathematical theory with spatial homogeneity} \label{Sec:4}

As a follow-up to Section 3, we consider a review and critical analysis for the study of systems in the spatially homogeneous case, where the microscopic variable is only the activity. Thus, the dependent variable describing the overall state of the system is, for each functional subsystem, the time-dependent one-particle distribution function over the activity.

The presentation follows a style similar to Section 3. It is divided into three subsections. First, we define the conceptual framework that leads to the derivation of a differential system considered to describe the dynamics of the system, which is essentially based on the dynamics of collective learning followed by decisions. Then, the presentation deals with the derivation of the differential system, which also considers the interactions within networks. Finally, a critical review of the literature on applications and parallel methods is developed. In this subsection, some research perspectives are indicated.

A preliminary consideration is that, although the differential system does not include the transport term in space, it does include several new features of the dynamics, such as proliferative and destructive interactions, corresponding to mutation-competition-evolution and transitions across functional subsystems. This rich framework explains why a wide variety  of applications have been developed in the life and applied sciences.

\subsection{New conceptual frameworks}\label{Sec:4.1}

This section provides the conceptual framework that leads to the derivation of mathematical structures.  We do start with the current state of the art and then, propose some developments. As in Section 3, we consider multiple and nonlinear interactions. In fact, the assumption of binary interactions, which should be  acceptable  for a rarefied gas, should be considered as a very special case in the case of living systems. This structure transfers the dynamics at the low individual based scale to the collective motion described by the kinetic theory of active particles.

We consider a system consisting of a number \textit{n} of FSs, denoted by the index $i= 1, \ldots, n$. First, we will consider systems with  a-particles whose activity is a scalar, while  the case of vector activity will be treated later. The state in each FS is described by the distribution function
\begin{equation}\label{FS}
f_i = f_i(t, u)\, : [0,T] \times D_u \, \longrightarrow \mathbb{R}_+, \hskip.5cm i= 1, \ldots, n,
\end{equation}
where  $D_u$ denotes the domain of the variables $u$.

\vskip.2cm \noindent \textbf{Remark 4.1.} \textit{The activity variable, in most cases, cannot be referred to well defined measure units, but it is possible observing their minimal and maximal values $u_m$ and $u_M$, so that u can be divided by $|u_M - u_m|$, so that if the activity is a positive quantities and is divided by $|u_M - u_m|$ then, one has $D_u = [0,1]$. If $u$ can attain negative and positive values symmetric with respect to $u=0$, the same considerations lead to  $D_u = [-1,1]$.}

\vskip.2cm The following particles are statistically involved in the dynamics:

\vskip.2cm \noindent \textbf{Test particle} of the $i$-th functional subsystem with microscopic state, at time $t$, delivered by the variable $u$, whose distribution function is $f_i = f_i(t, u)$. The test particle is assumed to be representative of the whole system.

\vskip.2cm \noindent \textbf{Field particles} of the $k$-th functional subsystem with microscopic state, at time $t$,
defined by the variable $u^*$, whose distribution function is $f_k = f_k(t, u^*)$. By field particles we denote the full set of particles.

\vskip.2cm \noindent \textbf{Candidate particles} of the $h$-th functional subsystem,  with microscopic state, at time $t$, defined by the variable $u_*$, whose distribution function is $f_h = f_h(t, u_*)$. Candidate particles are the field particles that end up in the test particle state after interaction.
\vskip.2cm

Interactions can induce a complex dynamics, which can generate:

\vskip.2cm \noindent (i) A modification in the microscopic state, i.e. the activity, which can undergo a progression/regression dynamics.

\vskip.2cm \noindent (ii) Proliferative and/or destructive events, this dynamic can also be related to a birth dynamic followed by mutations and Darwinian type selection.

\vskip.2cm \noindent (iii) Transition across functional subsystems. In fact, functional  subsystems can aggregate individuals with close affinity properties.

\vskip.2cm In general, according to the theory developed in~\cite{[BBGO17]}, the dynamic depends on the  microscopic state and distribution function of the interacting entities. Let us now provide some details on the above steps. We will reduce the overlap with the content of Section 3 as much as possible. However, some repetitions cannot be avoided in order to demonstrate the specificity of the dynamics under consideration.

\begin{enumerate}
\vskip.2cm  \item    \textit{Interaction Dynamics:} The modeling approach \textit{can be developed along two sequential steps}. First, the a-particles learn about the state of the other particles, acquiring (learning) the cues from external actions. Then they decide how to modify their activity, either by changing it, by promoting a proliferative/destructive action, or by moving to another FS.

\vskip.2cm  \item  \textit{Interactions and heterogeneity:} Interactions are nonlinear and multiple. Their modeling should take into account the specific strategy expressed by the a-particles. Such a strategy is heterogeneously distributed among them. The \textit{sensitivity domain} $\Omega$ is the domain of the activity variable within which an a-particle can feel the presence of the other a-particles.

\vskip.2cm  \item \textit{Multiscale interactions:} Two types of interactions are considered: \textit{micro-micro interactions} which involve the activity variable of the interacting particles, and \textit{micro-macro interactions}  which involve the activity variable of the  a-particle and the whole system represented by low-order moments of the of the system.

\vskip.2cm  \item The \textit{decision making} is induced by a utility function that evolves in time and then  acts over different strategies, which can include  ``rational'', but also ``irrational'' behavior. The decision selects the output of interactions by using an \textit{utility function} that considers not only the current state of the system, but also an estimate of the perspective evolution of the whole system.

\vskip.2cm  \item \textit{Competition}: When one of the interacting a-particles increases its state by taking advantage of the other, which is forced to decrease its state. \textit{Consensus}: When a-particles exchange their status, i.e. a-particles with higher status decrease it, while the others with lower status increase it.

\vskip.2cm  \item\textit{Learning:} One of the two a-particles modifies the micro-state independently of the other. It learns by decreasing the distance between them.

\vskip.2cm  \item\textit{Hiding}: One of the two tries to increase the total distance from the other, while in turn, the other tries to decrease it, called \textit{chasing}.

\end{enumerate}

The overall approach is that mathematical models related to living systems should take into account that the rules of interaction  are not fixed in time, but evolve over time within the so-called \textit{Artificial World} of  Herbert A. Simon. We refer to his philosophy, which starts from a complexity analysis and leads to a deep analysis of the dynamics of interactions, see~\cite{[Simon1965],[Simon1976]} and the mathematical interpretation given in~\cite{[BE24]}.

\vskip.2cm \noindent  \textbf{Remark 4.2.} \textit{In general, mixed interaction rules, such as hiding-chasing and competitive-cooperative, have been used systematically in modeling social dynamics, see~\cite{[ABG16],[DLO17]}. Other applications have developed the idea that consensus versus dissent depends on the social distance of the interacting entities. The mathematical interpretation shows that the asymptotic trend leads to different distributions of wealth, small or large size of the middle class compared to populations with low and high levels of wealth. Analogous behaviors appear in the study of clustering in opinion dynamics.}

\vskip.2cm \noindent  \textbf{Remark 4.3.} \textit{The book~\cite{[BL12]} develops a mathematical approach to the study of social and economic systems using the mathematical tools of systems theory. The last part of the book brings to the attention of the authors some research perspectives, in particular some open problems, which the authors themselves point out by analogy with the famous Hilbert's problems~\cite{[Hilbert]}. One of the remarks is that most of the mathematical tools are based on consensus dynamics, whereas different types of interactions should be studied. Indeed, it is true that consensus problems have been studied extensively, see the review~\cite{[MT14]} and the bibliography therein. On the other hand, Remark 4.2 shows that the current state of the art already presents studies on this specific topic, see~\cite{[APZ19],[DL15]}, as well as studies focused on algorithms~\cite{[HSK20],[HSK21]}.}

\vskip.2cm
Both of the above remarks indicate that multiple strategies may be present in the heterogeneous behavior of interacting a-particles. An interesting contribution to this topic is given by the so-called \textit{Parrondo's Paradox}, which shows that two different strategies may lose if applied separately to pursue the same objective, but may gain if applied together in some, perhaps convex, combination. A sharp analysis of this topic is developed in~\cite{[WC24]}.

\subsection{Derivation of mathematical structures}\label{Sec:4.2}

The derivation of the mathematical structures follows the same methodology that has been used in the case of space-dependent dynamics. In the case of spatial homogeneity, the derivation needs to be more detailed, as the variety of applications is large and has suggested further developments of the structure. First, we consider the modeling of interactions, specifically the structures that can describe these dynamics. Then,we consider the derivation of the structure and some further technical developments. The presentation will consider both micro-micro and micro-macro interactions.

\vskip.5cm \noindent \textbf{Modeling interactions:} In the case of \textit{micro-micro interactions}, a-particles interact with a-particles within an interaction domain $\Omega_i$, where the subscript corresponds to the FS.  The size of the interaction domain $\Omega_i$ depends on the amount of information that can be received by an a-particle, i.e.  it depends on the distribution function. However, unlike to the dynamics studied in Section 3, the activity is the only microscopic variable and a formalization of topological interactions is not currently available. Therefore, most applications take $\Omega_i$ coincident with the domain $D_u$ of the activity variable. Some authors have suggested using of a weight function to model the intensity of the interactions in $D_u$.
In the case of \textit{micro-macro interactions}, a-particles interact with each FS viewed as a whole. Then, with the representation of each functional subsystem by low order moments of the distribution function.

The modeling of the interactions is obtained by the following quantities:

\vskip.2cm \noindent \textit{Interaction rates:}  The frequency and the intensity of micro-micro interactions between a candidate $h$-particle  and a field $k$-particle with states $u_*$ and  $u^*$, respectively, is modeled by $\eta_{hk}[f_h,f_k](u_*, u^*)$. The frequency and intensity of interactions between a candidate $h$-particle with state $u_*$ and the $k$ FS  with $\EE_k$ is described by $\mu_{hk}[f_h, f_k](u_*, \EE_k)$.

\vskip.2cm \noindent \textit{Micro-micro interactions:} $\cA_{hk}^i[f_h,f_k](u_* \rightarrow u|u_*, u^*)$ models  the probability density that a  candidate $h$-particle, with state $u_*$, ends up into the state of the test particle of the  $i$-th FS after an interaction  with a field $k$-particle with state $u^*$.

\vskip.2cm \noindent \textit{Micro-macro interactions:} $\cM_{hk}^i[f_h,f_k](u_* \rightarrow u|u_*, \EE_k)$ models  the probability density that a  candidate $h$ a-particle, with state $u_*$, ends up into the state of the test particle of the  $i$-th FS after an interaction  with the $k$ FS.

\vskip.2cm \noindent \textit{Loss terms:}  $\cL_i $ and $\cL_i^M $ model the loss of number of particles corresponding to micro-micro and micro-macro interactions due to conservative interactions.

\vskip.2cm \noindent\textit{Proliferative and destructive terms:} We consider only micro-micro interactions under the assumption that both dynamics occur within the same FS of the test particles due to interactions with the field particles. In detail, the dynamics refer to a test $i$-particle interacting  with a field $k$-particle with state $u^*$. Then  $\cP_{ik}[f_i,f_k](u, u^*)$ and $\cD_{ik} [f_i,f_k](u, u^*)$ model the proliferative and the destructive events, respectively.

\vskip.5cm \noindent \textbf{From interactions to mathematical structures:}  The derivation of mathematical structures, which are considered to provide the conceptual framework for modeling, is obtained by the numerical balance of a-particles within an elementary volume of the space of microscopic states of the active particles $[u, u + du]$. The balance of particles in the elementary volume of the space of microscopic states leads to the said derivation. More precisely, the rate of change of the number of active particles is equal to the input flux rate minus the output flux rate for both micro-micro and micro-macro interactions. The numerical balance is analogous to that in Section 3 and, as mentioned above, we consider the case of proliferation and destruction only in the $i$-FS. Then, the technical calculations yield the following:
\begin{eqnarray}\label{general}
&& \frac{\p}{\p t} f_i(t, u) = Q_i[\bbf](t,u) \nonumber \\[2mm]
&& \hskip1.6truecm = \bigg((\cG_i -  \cL_i)  + (\cG_i^M -  \cL_i^M)  + (\cP_i - \cD_i)\bigg)[\bbf](t,u),
  \end{eqnarray}
where $\bbf = \{ f_i\}$, and
\begin{eqnarray}\label{mm-gain}
&& (\cG_i - \cL_i) = \sum_{h,k=1}^n \,\int_{\Omega_i \times \Omega_i} \eta_{hk}[f_h,f_k](u_*,u^*) \mathcal{A}_{hk}^i[f_h,f_k]\left(u_* \to u|u_*,u^* \right)\nonumber\\[2mm]
&& \hskip3cm \times f_h(t,u_*)f_k(t, u^*)\,du_* du^* \nonumber\\[2mm]
&& \hskip1.5cm  - f_i(t,u) \sum_{k=1}^n \, \int_{\Omega_i}  \eta_{ik}[f_i,f_k](u,u^*)\, f_k(t,u^*)\,du^*,
\end{eqnarray}
\begin{eqnarray}\label{mM-gain}
&& (\cG^M_i -  \cL_i^M) = \sum_{h,k=1}^n \,\int_{\Omega_i} \mu_{hk}[f_h,f_k](u_*,\EE_k) \cM_{hk}^i[f_h,f_k]\left(u_* \to u|u_*,\EE_k \right)\nonumber\\[2mm]
&& \hskip3cm \times f_h(t,u_*)\EE_k(t)\,du_*\nonumber\\[2mm]
&& \hskip1.5cm  - f_i(t,u) \sum_{k=1}^n \,\int_{\Omega_i} \mu_{ik}[f_i,f_k]\,(u,\EE_k) f_k(t,u^*)\,du^*,
\end{eqnarray}
and
\begin{eqnarray}
&&(\cP_i - \cD_i) = \sum_{k=1}^n \,\int_{\Omega_i \times \Omega_i} \eta_{ik}[f_i,f_k](u_*,u^*)\mathcal{P}_{ik} [f_i,f_k](u_* \to u|u_*, u^*) \nonumber\\[1mm]
&&\hskip3cm \times  f_h(t,u_*)f_k(t,u^*)\,du_* du^*\nonumber\\[2mm]
&& \hskip1.5cm - f_i(t, u) \sum_{k=1}^n \,\int_{\Omega_i} \eta_{ik}[f_i,f_k](u, u^*)\, \mathcal{D}_{ik}[f_i,f_k]\,f_k(t, u^*)\, du^*.
\end{eqnarray}

If an external action $\vf_i[{\color{red}\bbf}](t,u)$ is applied to the activity variable, as mentioned in Section 3, the structure of the left-hand side term modifies as follows:
$$
\frac{\p}{\p t} f_i(t, u) + \frac{\p}{\p u}\big(\vf_i[{\color{red}\bbf}](t,u) f_i(t,u)\big)  = Q_i[{\color{red}\bbf}](t,u).
$$

For the sake of completeness, we will report on two technical generalizations that have been applied in the context of modeling specific systems. In particular, we consider structures that use \textit{discrete activity states}, corresponding to modeling the activity variable by a finite number of activities and structures that use of \textit{vector activity variables}.

The search for mathematical structures with discrete activity states is motivated by the assumption of a continuous distribution function, which requires that the number of a-particles is sufficiently large to justify the continuity assumptions. In addition, the behavioral variables are difficult to measure with sufficient precision. Therefore, it is appropriate and useful to discretize the activity into a finite number of intervals.
This idea was proposed in~\cite{[BLS00]}, see Chapter 10, and applied by various authors as reviewed in~\cite{[ABG16],[DLO17]}. Analytical problems, including the development of the Hilbert problem, are treated in the book~\cite{[ABDL03]}.

New ideas to derive discrete models have been developed to model vehicular traffic, see~\cite{[CDF07]}, where the authors proposed the interesting idea that the size of the discretization should depend on the local density and be variable in time. In general, the motivation for using models with discrete values of the activity variable is compelling, based on the concept that the discrete representation of activity is appropriate for a variable whose measurements are more reliable when referred to intervals rather than a single value.

The papers cited above consider models with a constant number of particles, but the generalization to include proliferative and destructive interactions is immediate, as already dealt with in~\cite{[CDF07]} and further developed in~\cite{[BBD21]}.

Mathematical structures  can be derived from the above mentioned discretization, resulting in the following dependent variables:
\begin{equation}\label{dep-variabe}
f_{ij} = f_{ij}(t), \hskip.5cm \hbox{with} \hskip.5cm i= 1, \ldots n, \hskip.5cm j= 1, \ldots m, \hskip.5cm  \bbf = \{f_{ij}\},
\end{equation}
where the subscripts $i,j$ correspond to the functional subsystem and the discrete value of the activity. The following calculations:
\begin{equation}
n_i = n_i(t) = \sum_{j=1}^m f_{ij}(t) \hskip.5cm \hbox{and} \hskip.5cm  N= N(t) = \sum_{i=1}^n  \sum_{j=1}^m f_{ij}(t)
\end{equation}
denote the local number density of each FS and for the whole population respectively.

If we denote by $ij$-particle the a-particle with state $j$ in the i-th FS, then the interactions are modeled by the following quantities:

\vskip.2cm \noindent $\cA_{pq}^{hk}[\bbf](p \to i, q \to j)$ models the transition from the p-FS to the i-FS and from the $q$-state to the $j$ state. These dynamics refer to the micro-micro interactions between a $pq$-particle and a $hk$-particle, which interacting at the rate $\eta_{pq}^{hk}$.

\vskip.2cm \noindent $\cM_{pq}^{k}[\bbf](p \to i, q \to j)$ models the  analogous transition for the micro-macro interactions between a $pq$-particle and a and the $k$-FS, which interacting at the rate $\mu_{pq}^{k}$.

\vskip.2cm The derivation of the structure is developed under the assumption that proliferative and destructive events occur  in the same FS  and in the state of the test particle.

\vskip.2cm \noindent $\mathcal{P}_{ij}^{hk}[\bbf]$ models the proliferation of a ij-particle due to micro-micro interactions with  hk-particles.
 The  interaction rate is $\eta_{ij}^{hk}$.

\vskip.2cm \noindent ${D}_{ij}^{hk}[\bbf]$ models the  analogous destructive event, which occurs with  $\eta_{ij}^{hk}$.

\vskip.2cm
Then, technical calculations analogous to those used for the continuous distribution lead to the following general structure:
\begin{eqnarray}\label{discrete}
&&\frac{d}{dt} f_{ij} = \sum_{p,q = 1}^n \sum_{h,k = 1}^m  \eta_{pq}^{hk}[\bbf] \cA_{pq}^{hk}[\bbf](p \to i, q \to j) f_{pq} f_{hk}\nonumber \\[2mm]
&& \hskip2cm  - f_{ij}\sum_{h=1}^n \sum_{k = 1}^m \eta_{ij}^{hk}[\bbf] f_{hk} \nonumber \\[2mm]
&& \hskip1cm + \sum_{p,q = 1}^n \sum_{k = 1}^{m}  \mu_{pq}^k[\bbf] \cM_{pq}^{k}[\bbf](p \to i, q \to j) f_{pq} \EE_{k}\nonumber \\[2mm]
&& \hskip2cm  - f_{ij}\sum_{h=1}^n \sum_{k = 1}^m \mu_{ij}^k[\bbf] \EE_{k}\nonumber \\[2mm]
&& \hskip1cm  + \sum_{h = 1}^n \sum_{k=1}^n \ \eta_{ij}^{hk}[\bbf]\,\mathcal{P}_{ij}^{hk} [\bbf] f_{ij}\, f_{hk}\nonumber \\[2mm]
&& \hskip2cm  - f_{ij} \sum_{k=1}^n \, \eta_{ij}^{hk}[\bbf]\, \mathcal{D}_{ij}^{hk}[\bbf]f_{hk},
\end{eqnarray}
where all dependent variables $f_{ij}$ and related moments depend on time.

Specific calculations can generate either simpler structures or more general ones. As an example,
we can consider the case micro-micro interactions in absence both of  transitions across FSs  proliferative and destructive interactions. The structure is as follows.
\begin{eqnarray}\label{simple}
&&\frac{d}{dt} f_{ij} = \sum_{h = 1}^{n} \sum_{q,k = 1}^{m}  \eta_{iq}^{hk}[\bbf] \cA_{iq}^{hk}[\bbf](q \to j) f_{iq} f_{hk}\nonumber \\[2mm]
&& \hskip2cm - f_{ij}\sum_{h= 1}^n \sum_{k = 1}^m \eta_{ij}^{hk}[\bbf] f_{hk}.
 \end{eqnarray}

Occasionally, the analogy with the so-called \textit{Discrete Boltzmann Equation} is mentioned, see~\cite{[Gatignol],[MP1991],[Plat-Illner]}. The analogy is correct in the case of discrete velocity as in~\cite{[CDF07]} for vehicular traffic, and in~\cite{[BBK13]} for crowd dynamics, but it is rather artificial in the case of discrete activity.

The case of \textit{vector activity variables} can be treated, on a formal level, simply by replacing the scalar $u$ by the vector $\bu = \{u_1, \ldots, u_\kappa\}$. However, this general approach can cause difficulties in the modeling of interactions.

The hint suggested in~\cite{[BBD21]} is to model a hierarchy
$$
u_1 \to u_2,  \hskip1cm u_1, u_2 \to u_3,  \hskip1cm u_1, u_2, u_3 \to u_4, \ldots
$$
This assumption may allow a factorization of the transition density. The drawback is that it is not a general strategy. Then it can  only be used in specific case studies.

\subsection{From theory to models - A narrative description}\label{Sec:4.3}

We present an overview of applications based on the mathematical tools derived in the previous subsections. The presentation concerns applications related to a well-defined research area. A narrative presentation is adopted to highlight how each application has developed novel mathematical methods and tools. Special attention is given to the use of multiple nonlinear interactions and multiple strategies. Furthermore, we pay attention to the role of collective learning dynamics and of the subsequent decision making that induces the collective motion.

The presentation developed in the following paragraphs is limited to a review and critical analysis of applications with some hints on developments in each specific research area.  Further perspectives are discussed in the subsequent  sections of this paper. The presentation does not claim to be exhaustive. Rather, we focus on articles that mark advances in mathematical tools. The applications have been focused in different fields, such as the immune competition, epidemics and economics. These applications have already been reviewed in surveys, see~\cite{[BBD21],[DLO17]} and in the short book~\cite{[BBST24]}. Therefore, we will not repeat the review, but, as an alternative, we will focus on two specific topics, i.e., the hints that applications have given to the development of the theory and research papers authored by an interdisciplinary team.

The first topic is presented with focus on the study of the \textit{immune competition between tumor and immune cells.} The research activity developed in this specific area was promoted by the seminal paper~\cite{[BF94]}, where two FSs were taken into account, i.e. tumor and immune cells. The model accounts for  the ability of immune cells to detect the presence of the tumor. This learning dynamic activates immune cells, which progress from the state of innate immunity to that of active immunity and proliferate up to a possible plafond corresponding to the highest activation. On the other hand, tumor cells progress their activation up to metastatic competence and proliferate by invading the surrounding tissue.

The most significant results are reported and studied in the book~\cite{[BD06]}, where the competition depends on a parameter that measures the ration between the activation ability of the immune cells and the progression ability of the tumor. This parameter separates the two different asymptotic trends corresponding to the gain of the immune system with depletion of the tumor system or the opposite corresponding to the gain of the tumor cells.
This excellent book also has the merit of having developed analytical studies on the solution of the differential systems describing  the above dynamics and of having developed the methodology to derive macroscopic tissue models from the underlying description provided by the kinetic theory models at the microscopic scale.

A number of  papers have been proposed on this line of research, focusing on the modeling topic, for example see~\cite{[FGP99],[LoSchiavo]}, or on the micro-macro derivation, see~\cite{[BBC16],[BD11]}. For the improvement that takes into account the modeling of tumor cells to mutate and the related ability of immune cells to learn the presence of the tumor, see~\cite{[BDK13],[CGS18]}. The biological theory related to cancer phenomena is reported in~\cite{[HW11]} and~\cite{[WEI07]}, see also~\cite{[HART01]}, while a general theory of immune competition is available in the book~\cite{[MF19]}.

However, further significant improvements in the kinetic theory approach did not follow in the above line of research, except for technical modifications that are not reported here. Although the success in vaccination studies by the Nobel laureates, in particular James P. Allison and Tasuku Honjo \textit{for their discovery of cancer therapy by inhibition of negative immune regulation}, see~\cite{[AH18]}.

The success of vaccination trials should motivate further mathematical studies to model the complex dynamics of immune system activation. Recent studies on modeling collective learning dynamics, see~\cite{[BDG16],[BDG16A]}, and nonlinear interactions, see~\cite{[BBD21]}, are supposed to contribute to further developments in cancer modeling. What mathematicians might learn from the above research is that a multiscale approach is necessary and that modeling should include dynamics at the lower molecular (genes) scale. Mathematicians have started to move towards this perspective, which will definitely become the key interpretation of multiscale vision, see~\cite{[CJL19],[CL06]}.

Further clues were provided by mathematical studies of the SARS-CoV-2 pandemic, which went far beyond the deterministic population dynamics approach to epidemics. The study of the in-host dynamics, which has been studied by active particle methods, has provided a detailed description of the pathology including the modeling of the level of pathology and hence the level of infectivity in a heterogeneous population, see~\cite{[BBC20],[BBO22],[BK24]}, see~\cite{[B3EPT]} for a review of the literature and the further developments of the kinetic theory of active particles. It is important to mention that the kinetic equations have been coupled to the network of transportation means to obtain the modeling of the spread of epidemics in the territory, see\cite{[BP21],[BDP20]}, an important reference on this topic is~\cite{[PCVV15]} that anticipates the problems of SARS-CoV-2 pandemics. In this framework, developing mathematical tools to measure heterogeneity is an important topic, see~\cite{[Tos25]}.

Another area of research, where active particle methods have been applied (along the lines of the above articles) is the modeling of the study of social dynamics.  Several authors have contributed to this research topic using models with linear interactions. We will skip this literature as well, since it is reported in previous reviews such as~\cite{[DLO17]}. On the other hand, we will focus on models with nonlinear interactions and/or models where multiple interactions occur. In fact, these topics were mentioned in the conclusions of the book~\cite{[BL12]}, where the authors proposed some key issues. One of them is the idea of going beyond the use of consensus dynamics to search for different models of interactions. This hint was timely, as some papers were already moving in this direction. In particular, the study of how a welfare policy can influence or oppose governments was proposed in~\cite{[BHT13]}.

An important feature of the modeling of the interactions was the use of the role of the mean value of the variable modeling the said support/opposition. In fact, such mean value can act as an attractor for the heterogeneous political tendency of the population. Simulations in~\cite{[BHT13]} have shown that this attraction can bring the whole population to a radicalization of the opposition, interpreted as a black swan~\cite{[Taleb07]}. Indeed, this concept has already been proposed in~\cite{[BS12]} as an attractor of the interacting individuals in swarms.

These seminal papers fostered a new trend in research activity, which is documented in some subsequent papers. For example, for the interplay between altruism and selfishness see~\cite{[DL14]}, while the role of networks in enhancing consensus has been studied in~\cite{[DL15]}. Additional studies focusing on the interactions between economic development and political competition have been developed in the context of escaping the \textit{blocking} trap, see~\cite{[DKLM17]}, while the dynamics of liquidity profiles on interbank networks have been modeled in~\cite{[DLM21]}. Further references are given in the review paper~\cite{[DLO17]}. The book by Giovanni Dosi, see~\cite{[DOSI2023]}, reviews and critically analyzes the study of idiosyncratic learning and technological development through mathematical models, see also~\cite{[DPV17]}. In particular, kinetic theory models have been developed in~\cite{[BDKV20]}.
*
Parallel studies of social dynamics and economics have been developed by means of agents methods, see~\cite{[BMST15],[DR2019],[Galam],[MMT19]}, deterministic particles, see~\cite{[FR21]}, dynamics over networks, see~\cite{[KNO13],[KNO14]}.  See also~\cite{[KTS20]}, where the modeling was developed by applying of the mathematical theory of  behavioral swarms, which will be reviewed in the following.

The above literature, generally developed through interdisciplinary collaborations, refers to behavioral economics, see~\cite{[Thaler16],[TS16]}, but also to political interpretations, see~\cite{[Schump47]} and~\cite{[AR06]}. Historically, we know that interdisciplinary collaborations have led to new mathematical theories, such as the so-called \textit{Nash equilibria}, see \cite{[NASH51],[NASH96]}, and the the theory of differential games, see~\cite{[F71],[F72],[I65],[LL06i],[LL06ii],[LL07]} for the theory, and, for the specific applications, see~\cite{[Achdou22],[Bartucci25]}, that are selected among others. Further parallel frameworks to be considered are the \textit{Fokker-Plank-Boltzmann}, see~\cite{[PT13]}, and the \textit{behavioral swarms theory}, see~\cite{[BHLY24]}. The mathematical theory of active particles should be related to the above mentioned frameworks, which will be reviewed in the next section.

The following two remarks focus on challenging  multiscale and interdisciplinary problems, see~\cite{[B3EPT]}. Figure 4 shows some aspect of these dynamics.
\vskip.2cm \noindent  \textbf{Remark 4.4.} \textit{Application to modeling immune competition against a proliferating virus. The in-host model, see \cite{[BBO22],[BK24]}, has been developed on the scale of cells and viral particles. This has led to a description that goes beyond that of population dynamics, as the model also describes the progression or regression of pathology at the level of individuals. These move in the territory with a spatial dynamics that also depend on transport networks, see~\cite{[BDZ25]}.}

\vskip.2cm \noindent  \textbf{Remark 4.5.} \textit{Interactions with the external world have been outlined in Section 3 following Remark 3.3, and  we can reiterate that dealing with this topic will require a deep knowledge of the physics of each system under consideration and a multiscale interpretation, which should, arguably, lead to new mathematical structures.}

\vskip.2cm
Some aspects of the interaction with the outside world are shown in Figure 4, which illustrates some of the multiscale and interdisciplinary features of the system. In particular, the figure shows that the national service is concerned with the health of infected people and vaccination programs. In this activity, not only the number of infected people is important, but also the knowledge of their pathological state plays a role in planning the hospitalization program. This specific information can contribute to the planning of vaccination programs. These generate another dynamic concerning the social interactions between the crisis managers of the epidemics and the citizens, who can either accept or oppose the vaccination programs. A more detailed description is proposed in~\cite{[B3EPT]}.

\begin{figure}
\tikzstyle{decision} = [diamond, draw, fill=blue!20,
    text width=6.5em, text badly centered, node distance=3cm, inner sep=0pt]
\tikzstyle{block1} = [rectangle, draw, fill=blue!20,
    text width=6.5em, text centered, rounded corners, minimum height=4em]
    \tikzstyle{block3} = [rectangle, draw, fill=green!20,
    text width=6.5em, text centered, rounded corners, minimum height=4em]
    \tikzstyle{block2} = [rectangle, draw, fill=red!60,
    text width=6.5em, text centered, rounded corners, minimum height=4em]
     \tikzstyle{block4} = [rectangle, draw, fill=black!60,
    text width=6.5em, text centered, rounded corners, minimum height=4em]
    \tikzstyle{block6} = [rectangle, draw, fill=yellow!60,
    text width=6.5em, text centered, rounded corners, minimum height=4em]
     \tikzstyle{block8} = [rectangle, draw, fill=violet!60,
    text width=6.5em, text centered, rounded corners, minimum height=4em]
    \tikzstyle{block9} = [rectangle, draw, fill=white!60,
    text width=6.5em, text centered, rounded corners, minimum height=4em]
    \tikzstyle{block10} = [rectangle, draw, fill=white!60,
    text width=6.5em, text centered, rounded corners, minimum height=4em]
\tikzstyle{line} = [draw, -latex']
\tikzstyle{cloud} = [draw, ellipse,fill=red!20, node distance=3cm, minimum height=2em]
\begin{center}\scalebox{0.85}{
\begin{tikzpicture}[node distance = 3.2cm, auto]
    \node [block2] (Inf) {Low\\ Infection\\ home care};
    \node [block1, left of=Inf] (Hea) {Healthy \\ population};
    \node [block2, right of=Inf] (Hos) {Middle\\Infection\\ Hospital care};
    \node [block6, above of=Hea] (Vacc) {Vaccination \\ Programs};
    \node [block6, below of=Hea] (Long) {Contagion \\ by contact \\and mobility};
    \node [block3, above of=Hos] (Rec) {Social \\ Health \\ Service};
    \node [block8, right of=Hos] (I-Care) {High Level \\ Infection \\ Intensive care};
    \draw [line width=.5mm, red, ->] ([yshift=0.3cm] Hea.east) --  ([yshift=0.3cm] Inf.west);
    \draw [line width=.5mm, blue, ->] ([yshift= - 0.3cm] Inf.west) --  ([yshift= - 0.3cm] Hea.east);
    \draw [line width=.5mm, red, ->] ([yshift=0.3cm] Inf.east) --  ([yshift=0.3cm] Hos.west);
    \draw [line width=.5mm, blue, ->] ([yshift= - 0.3cm] Hos.west) --  ([yshift= - 0.3cm] Inf.east);
    \draw [line width=.5mm, red, ->] ([yshift=0.3cm] Hos.east) --  ([yshift=0.3cm] I-Care.west);
    \draw [line width=.5mm, blue, ->] ([yshift= - 0.3cm] I-Care.west) --  ([yshift= - 0.3cm] Hos.east);
    \draw [line width=.5mm, blue, ->] (Rec.south) --  (Inf.north);
    \draw [line width=.5mm, blue, ->] (Rec.south) --  (Hos.north);
    \draw [line width=.5mm, blue, ->] (Rec.south) --  (I-Care.north);
    \draw [line width=.5mm, blue, ->] (Rec.west) --   (Vacc.east);
    \draw [line width=.5mm, blue, ->] (Vacc.south) -- (Hea.north);
    \draw [line width=.5mm, red, ->]  (Long.north) -- (Hea.south);
   \end{tikzpicture}}
\end{center}
\begin{center}
\caption{Multiscale and interdisciplinary aspects of virus epidemics}\label{fig4}
\end{center}
\end{figure}
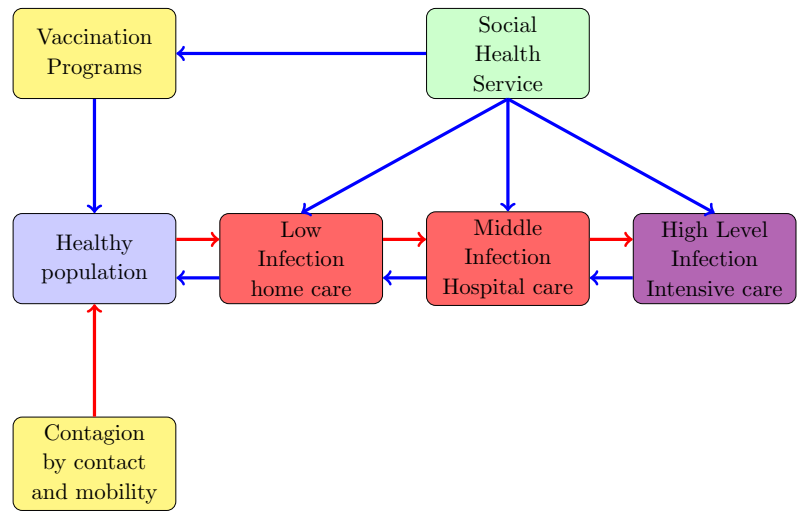

\section{Further parallel mathematical tools} \label{Sec:5}

The mathematical theory reviewed and critically analyzed in the previous sections refers to the kinetic theory of active particles, which is not the only method for dealing with collective dynamics. It is a method that has been developed and put into a general framework inspired by the seminal work of Ilia Prigogine. The key step to active particle methods is given in the paper~\cite{[BF94]} and the book~\cite{[BLS00]}. However, it should be noted that alternative, somewhat different approaches have been developed by other authors. Although they are quite different from the KTAP approach, we think that a brief account of them may be useful to provide a general vision of a research topic that has attracted applied mathematicians to intense studies in this century.

Specifically, we provide a concise overview of three methods that share a focus on the dynamics of living, i.e., complex, systems, although they are technically quite different. In detail, we consider the \textit{theory of mean field games} introduced in~\cite{[LL07]} and further developed even in recent work~\cite{[Achdou22],[Bartucci25]}. Then, as a second topic, we consider the \textit{Fokker-Plank-Boltzmann} method, the main results on this topic are reported in the book~\cite{[PT13]}. Applications have been developed in various fields such as social dynamics, opinion formation, and behavioral economics. The third topic is \textit{behavioral swarm theory}, which was introduced in~\cite{[BHO20]}, developed in~\cite{[BHLY24]}, and applied in economics~\cite{[KTS20]} and crowd dynamics~\cite{[KLMY25]}.

The selection of these topics is based on the authors' research activities and interests. We do not claim to present a complete vision, nor do we place the three methods in competition with each other. In fact, the mathematical theories mentioned above have been developed under different motivations and philosophies, which can be mutually integrated among them. Indeed, a challenging perspective is the search for a unified philosophy and thus for a unified mathematical theory. We focus on the essential bibliography, leaving out minor technical modifications of the theories in papers that claim to have proposed a new theory.

\subsection{Mean field games}\label{Sec:5.1}

Mean Field Games (MFGs) deal with the study of strategic decision making in differential games with a large number of indistinguishable, rational, and heterogeneous {\it players} (corresponding to {\it a-particles} in our context), which lies at the intersection of game theory with stochastic analysis and control theory \cite{[LL07]}. These methods approximate Nash equilibria for differential games with symmetric interactions between players. In contrast to classical game theory, MFGs model the interaction of a representative player with the collective behavior of the other players, rather than collecting detailed state information about all other players, which is unrealistic for dynamic games with a large number of players.

The use of the term \textit{mean field} is inspired by the classical mean field theory of Vlasov in the study of statistical physics \cite{[Vla68]}. The primary focus of Vlasov was to describe plasmas in which the ions never collide, but instead have long-range interactions that cannot be properly described by Boltzmann's equation. This is the same setting for most differential games or systems of a-particles, where the players or a-particles do not collide, but interact with others through non-local long-range effects that can be similarly described by a force field related to the distribution of the system, then the evolution of the distribution function can be governed by a Vlasov-type equation. In this framework, each player (a particle) acts according to its minimization or maximization problem, taking into account the decisions of other players, and since the population of players is large, one can assume that the number of players goes to infinity and that a representative player exists, then the behavior of each individual is governed by the rest of the players via the mean field.

The development of mean field game theory began with a series of papers by Lasry and Lions \cite{[LL06i],[LL06ii],[LL07]} around 2005. Related ideas were developed independently and at about the same time by Caines, Huang, and Malham\'e \cite{[HCM1],[HCM2]}, under the name of the Nash certainty equivalence principle. Note that in traditional game theory, the subject of study is usually a two-player game in discrete time space, and the results can be extended by induction to more complex situations. However, for differential games or stochastic differential games with continuous states in continuous time, this strategy cannot be used because of the complexity generated by the dynamic interactions. On the other hand, mean-field games can handle large numbers of players through the mean representative player and simultaneously describe the complex dynamics in differential games or stochastic differential games \cite{[Del20]}, see also~\cite{[Bartucci19],[CD18]}.

Mean field game theory also lies at the intersection of mathematics and complex systems, it has been successfully developed in several different areas of mathematics such as probability theory, stochastic control, partial differential equations, calculus of variations, optimal transport, and requires innovative numerical solutions, with links and applications to other scientific fields such as biology, social sciences, economics, and financial markets, which model a wide variety of phenomena in complex systems with large numbers of particles, see~\cite{[Achdou22],[BCL21],[BBLL24],[Bartucci25],[Del20],[KTS20]}.

\subsection{Fokker-Plank-Boltzmann mathematical tools}\label{Sec:5.2}

An additional topic to be outlined is a parallel method, let's call it the Fokker-Plank-Boltzmann method, by which Boltzmann-type kinetic models are first constructed for interacting multi-agent systems, and then Fokker-Plank-type equations are derived by quasi-invariant asymptotic procedures to capture the long-term behavior of the systems \cite{[CPT05],[FPTT17],[PT13],[T06],[TSB22]}. In these references, the authors use the term multiagent systems for interacting agents, which are the counterparts of the systems of active particles, and the interacting agents are the counterparts of the active particles.

Instead of considering the interaction rate and the transition probability density of the microstates, the Fokker-Plank-Boltzmann method follows the kinetic description of the interactions in the spirit of Boltzmann's ideas. More precisely, based on prescribed microscopic laws of {\it binary} interactions, Boltzmann-type equations are derived to describe the evolution of the one-particle density of some selected agent property, such as one's wealth \cite{[CPT05]} or opinion \cite{[T06]}, in a system of interacting particles or agents.

As an example, let $v\in \mathcal{I} $ be the value of the selected trait of the agents. Here the trait denotes a distinguishing quality or characteristic of the agents, such as wealth, opinion, knowledge, belief, or others. The domain $\mathcal{I} \subset \R$ can be a fixed interval or the positive half-line, depending on the application. Let $(v,w)$ be the values of the selected attribute of the pair of agents before their interaction, their post-interaction values $(v^*, w^*)$ can be given by the interaction rule \cite{[PT13]}
\begin{equation}\label{fp:int}
\begin{cases}
v^* =v+ P(v)(w-v)+Q(v)\eta,\\[2mm]
w^* =w+P(w)(v-w)+Q(w) \tilde{\eta},
\end{cases}
\end{equation}
where the non-negative functions $P$ and $Q$ contain the details of the interactions that measure the intensities of the variation of the trait due to the presence of the other agent, $\eta$ and $\tilde{\eta}$ denote the variations of the trait due to random effects that are assumed to be independent and identically distributed with zero mean and bounded variance.

Note that the functions $P$ and $Q$ and the random parameters are subject to conditions that keep the post-interaction values $(v^*, w^*)$ in the same domain $\mathcal{I}$ during the interactions. Then the time evolution of the distribution of the selected feature $v$ resulting from interactions of type \eqref{fp:int} between individuals can be described by collision-like kinetic models. That is, let $f = f(t,v)$ be the density of agents at time $t > 0$ with trait $v \in \mathcal{I} $, then it satisfies a weak form of a nonlinear Boltzmann-like equation, i.e. for all smooth functions $\varphi(v)$ (the observable quantities),
\begin{eqnarray}\label{fp:weak}
&&{d\over dt} \int_{\mathcal{I} } f(t,v)  \varphi(v) dv \nonumber\\[2mm]
&& \hskip1cm = {1 \over 2}\langle  \int _{\mathcal{I} \times \mathcal{I}}(\varphi(v^*) + \varphi(w^*) - \varphi(v) - \varphi(w) )
f(t,v) f(t,w) dv dw  \rangle,
\end{eqnarray}
where $\langle \cdot \rangle$ represents the mathematical expectation taking into account the presence of the random parameters $\eta$ and $\tilde{\eta}$ in (\ref{fp:int}). This kinetic equation means that the time variation of the distribution of characteristic $v$ is measured by the change in $v$ resulting from interactions of type (\ref{fp:int}).

Similar kinetic equations of Boltzmann type can also be obtained by assuming that the change of the trait $v \in \mathcal{I} $ is the consequence of interactions with a fixed environment. In this case, the variation of the individual trait in each single microscopic interaction is the result of the interaction rule \cite{[FPTT17]}
\begin{equation}\label{fp:int2}
v^* = v + P_E(v)z - P(v)v + Q(v)\eta,
\end{equation}
where the intensity of the non-negative variation $P_E(v)$ of the trait $v$ due to the presence of the environment is generally assumed to be different from $P(v)$, $z \in \mathcal{I} $ denotes the amount of trait absorbed by the agent from the environment, and $\eta$ denotes the random variations of the trait as before. Using the same notation as in \eqref{fp:weak}, in this case the density $f(t,v)$ satisfies, for all smooth functions $\varphi(v)$ (the observable quantities),
\begin{equation}\label{fp:weak2}
{d\over dt} \int _{\mathcal{I} } f(t,v)  \varphi(v) dv =
\langle  \int _{\mathcal{I} \times \mathcal{I} }
(\varphi(v^*) - \varphi(v)  )
f(t,v) {\mathcal{E}}(z)  dv dz  \rangle,
\end{equation}
where $\mathcal{E}(z), \ z \in  \mathcal{I}$ denotes the distribution of the environment.

Then next, the long-time behavior of the kinetic equations (\ref{fp:weak}) and (\ref{fp:weak2}) can be governed by Fokker-Planck-type equations \cite{[CPT05],[T06]}, which are derived by some well-consolidated asymptotic procedures which are reminiscent of the so-called grazing collision limit of the classical Boltzmann equation \cite{[V02]}, and the dissipative versions of Kac's caricature of a Maxwell gas, see \cite{[FPTT12],[PT04]}.

The analysis of the resulting Fokker-Planck-type equations can well explain many interesting phenomena in socio-economic applications \cite{[CPT05],[PT13],[T06],[Tos25]}. However, the basis of the Fokker-Planck-Boltzmann method relies more on binary interactions in the spirit of Boltzmann's ideas, which is quite different from the general mathematical structure for active particles as discussed in section \ref{Sec:3.2}. A more detailed consideration of this parallel Fokker-Plank-Boltzmann method is beyond the scope of the present report, and we refer to \cite{[FPTT17],[FPTT20],[PT13],[TTZ18],[ZBT21]} and references therein for more details on the analysis and applications of this approach.

\subsection{On the mathematical theory of behavioral swarms}\label{Sec:5.3}

Note that the above methods have been applied to the study of the dynamics of systems with a large number of interacting entities. The common feature of these methods is that the collective behavior of the whole system is described by probability distributions over the microscale state of the interacting entities, i.e., a-particles, players, or agents in different contexts. As an alternative approach, the \textit{behavioral swarm theory}, first introduced in~\cite{[BHO20]} and further developed in~\cite{[BHLY24]}, is useful for describing the dynamics of living systems with a finite number of interacting entities without applying the assumption that the number of interacting entities tends to infinity.

Modeling the dynamics of swarms with internal variables can be traced back to~\cite{[BS12]}, where the authors proposed introducing a {\it activity variable} related to the strategy expressed by individual entities in swarms. This variable was thought to model different emotional states that occur, for example, in predator-prey swarms, and a more general mathematical approach to behavioral swarms was established in \cite{[BHO20]}. This approach extends the classical Cucker-Smale theory \cite{[CS07]} to include the dynamics of internal variables that have the ability to interact with mechanical variables and thus affect the interaction rules. The modeling of interactions proposed in~\cite{[BS12]} has transferred into a mathematical framework the topological interactions conjectured in~\cite{[BCC08]}. Further interpretation of~\cite{[BS12]} have been developed  in~\cite{[ALBI24]}.

More specifically, in the mathematical theory of behavioral swarms, the state of interacting entities, also called a-particles, includes an internal variable -- activity -- in addition to the mechanical variables -- position and velocity. This additional activity variable has the ability to interact with the mechanical variables and thus affect the interaction rules. In turn, the mechanical variables can modify the dynamics of the activity variable. This approach is useful for describing the dynamics of living systems with a finite number of interacting entities.

Moreover, behavioral swarm theory extends the concept of \textit{swarm intelligence}, which is generally focused on the specific behaviors of animal swarm entities, to that of collective human behaviors. This extension can be seen as a natural step towards the methods of artificial intelligence \cite{[BHLY24]}. In fact, the philosophy of the proposed theory, although it follows the classical path, i.e.
$$
\mbox{\bf following} \hskip.2cm \longrightarrow  \hskip.2cm  \mbox{\bf interaction}  \hskip.2cm  \longrightarrow  \hskip.2cm  \mbox{\bf learning}  \hskip.2cm \longrightarrow  \hskip.2cm  \mbox{\bf decision making},
$$
takes into account the specificity of living systems. As a result, interactions are not simply based on \textit{consensus}, but may be somewhat different. Rather, it is decision making that evolves over time, taking into account the overall dynamics of the swarm.

As an application, the crowd dynamics with social interaction is considered in  \cite{[KLMY25]}, in which the dynamics of the activity variable can interact with that of the mechanical variables thus affecting the interaction rules in crowd dynamics with social interaction. This approach is useful for describing heterogeneous behavioral features in crowds, especially, the dynamics of pedestrians with evolving psychological states.

Further applications can also be addressed in different types of systems. For instance, the competitive in-host dynamics between invasive virus particles and of the immune system agents~\cite{[BBO22],[BK24]}, the dynamics of robots or unmanned vehicles in engineering sciences~\cite{[Shain2005]},  dynamical systems in economics and financial markets~\cite{[BCLY17],[BCK19],[KTS20],[WBZ23]}, in which different aspects to model causality principles have been studied. Actually most of the previous studies were developed under the assumption of spatial homogeneity. On the other hand  the general theory can consider more general inhomogeneous dynamics.

\section{A forward look to new research horizons}\label{Sec:6}

The mathematical literature reviewed in the previous sections provides a general description of the development and applications of active particle methods in the context of the collective dynamics of living entities. A critical analysis and some pointers to developments are also reported. Research perspectives are then brought to the attention of the reader. A further selection of these is proposed in the following subsections. We begin with a topic that we consider very important, namely the artificial world of Simon's philosophy. In fact, the mathematical interpretation of Simon's theory permeates the entire content of the following subsections.

Therefore, this section looks forward to research perspectives proposed by concepts and hints, while their development will of course include analytical and computational methods. The presentation will show how all topics are interrelated and tend towards a unified vision of a general mathematical theory of the collective dynamics of living systems. The goal is to understand to what extent the mathematical theory in~\cite{[BBGO17]} can be developed to include the dynamics of the artificial world and, ultimately, how this leads to the collective dynamics of living systems.

\subsection{Interactions in the artificial world of Herbert A. Simon}\label{Sec:6.1}

A challenging research perspective is the development of mathematical tools to describe the dynamics of evolutionary systems within the general framework of the philosophy of Herbert A.~Simon, see~\cite{[Simon1965],[Simon1976],[Simon2019]}. This study was initiated in~\cite{[BE24]} and is based on a new interpretation of interactions as a preliminary step in the pursuit of the above goal. This paper focuses on the dynamics of behavioral systems of evolutionary economics~\cite{[Egidi22],[EMS24]}. In particular, the theory suggests that, in the case of living systems, the search for universal interaction rules governing interactions should be conducted within a time-evolving environment, which he calls \textit{artificial world}. The events of the dynamics were called \textit{artifacts}. Some applications have been critically examined with a focus on self-organization in various types of aggregated systems, see~\cite{[EP92],[Felin-Foss],[Foss2012],[Ostrom57],[Ostrom98],[Ostrom2000]}.

It is worth emphasizing that the authors in~\cite{[BE24]} observe that some concepts of the complex connections between the real and sensitive worlds are present in the dissertation~\cite{[KANT]}, which marks the transition to the ``critical period'' leading to the dissertation, see also~\cite{[KANT-B]}. This dissertation begins with a definition of ``mundus'' and then defines the difference between ``\textit{mundus sensibilis}'' and ``\textit{mundus intellegibilis}''. An important feature of the \textit{mundus sensibilis} is the heterogeneity of the individual learning process.

This challenging perspective leads us to go far beyond the mathematical tools reviewed in the previous sections. Indeed, although we start from the kinetic theory of active particles~\cite{[BBGO17]}, the description of the interactions must be rethought, since the rules that guide them are not independent of the properties of the environment in which the dynamics develop (learning, decision making, development of strategies, etc.), but refer to the \textit{virtual world} that hosts them. In fact, the \textit{virtual world} follows its own dynamic rules depending on external actions, for example, due to actions of politics, as well as on internal dynamics of the different components of such world.

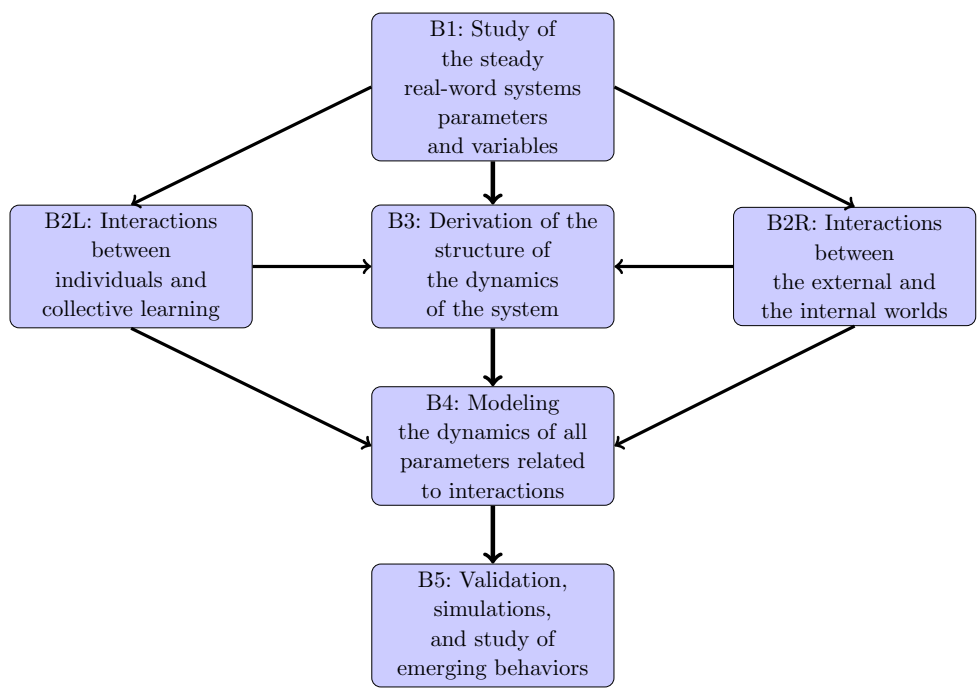
\begin{figure}[t!]
\begin{center}\scalebox{0.85}{
\begin{tikzpicture}[node distance = 2.8cm, auto]
 \node [block_long, below of=system] (complex) {B1: Study of  the steady \\real-word systems\\ parameters and variables};
 \node [block_long, below of=complex] (model-int) {B3: Derivation of the \\structure of the dynamics \\ of the system};
 \node [block_long, right of=model-int,xshift=2.8cm] (outer) {B2R: Interactions \\between\\ the external and \\the internal worlds};
 \node [block_long, left of=model-int,xshift=-2.8cm] (learning) {B2L: Interactions \\ between \\ individuals and \\ collective learning};
 \node [block_long, below of=model-int] (art-world) {B4: Modeling \\ the dynamics  of all\\parameters related \\ to interactions};
 \node [block_long, below of=art-world] (models){B5: Validation, simulations,\\ and study of \\emerging behaviors};
  \draw[line width=.7mm, ->] (complex.south) -- (model-int.north);
  \draw[line width=.5mm, ->] (complex.west) --  (learning.north);
  \draw[line width=.5mm, ->] (complex.east) --  (outer.north);
  \draw[line width=.5mm, ->]  (learning.south) -- (art-world.west);
  \draw[line width=.5mm, ->]  (outer.south) -- (art-world.east);
  \draw[line width=.7mm, ->]  (model-int.south) -- (art-world.north);
  \draw[line width=.5mm, ->] (outer.west) -- (model-int.east);
  \draw[line width=.5mm, ->] (learning.east) -- (model-int.west);
  \draw[line width=.7mm, ->] (art-world.south) -- (models.north);
\end{tikzpicture}
}
\end{center}
\caption{Strategy towards a mathematical theory and derivation of models}\label{fig5}
\label{rationale}
\end{figure}

The method is shown in the flowchart in Fig.~5, which can be interpreted as an evolution of the previous one in Fig.~1. The contents of the blocks are as follows:

\vskip.2cm \noindent \textbf{Block 1} focuses on the preliminary study of the system object to be modeled. Specifically to the systems discussed in Sections 3 and 4 and to the general philosophy in Section 2. The study identifies the variables and parameters needed to describe the system object of the modeling approach and the role they play in the dynamics. Both spatial dynamics and spatial homogeneity can be considered.

\vskip.2cm   \noindent \textbf{Blocks 2L and 2R} consider interactions at different scales, i.e. micro-scale interactions involving individual a-particles and micro-macro-scale interactions referring to each FS as a whole acting on individual a-particles. Block 2L refers to the interactions within the functional subsystem, while Block 2R refers to the modeling of the external actions on the a-particles in functional subsystems. The implementation of both descriptions contributes to the definition of the mathematical structure sought in Block 3.

\vskip.2cm \noindent \textbf{Block 3} essentially refers to the approach discussed in the previous sections, which
assumes that all a-particles interact with the particles in each FS according to rules that may depend on the specificity of each FS but are constant in time. The derivation of the mathematical structure should include all interactions associated with Blocks 2L and 2R. This structure is needed to capture as much as possible the complexity characteristics of living systems.

\vskip.2cm  \noindent   \textbf{Block 4} focuses on the most difficult part of the approach as it aims to derive a differential system to describe the dynamics of the parameters in relation to the overall evolution of the system. The derivation requires an interpretation of the role of the virtual world in determining the dynamics of the parameters, which characterize quantitatively the interactions.

\vskip.2cm  \noindent \textbf{Block 5} deals with the validation of the models. As mentioned in the previous sections, the derivation of specific models is obtained by inserting the mathematical description of the interaction into the structure mentioned in Block 3. The models have to be validated. This is not an easy task, considering that the study should be related to the dynamic behavior. This has already been addressed in~\cite{[BBGO17]}, where validation was related to the ability of the models to reproduce, at a qualitative level, the emerging behaviors in the collective motion of living systems.

\vskip.2cm
The topic covered in this subsection permeates the content of the following sections. In fact, it affects the entire process of learning and decision making, which is the rationale that guides the entire content of this paper. In particular, we expect that the next section will propose a critical analysis that gathers all the clues proposed in this paper towards the challenging goal of developing a mathematical theory of living systems.

\subsection{From collective learning to decision making}\label{Sec:6.2}

The previous sections have indicated that we look at a dynamics of the type:
\vskip.1cm
 \begin{center}
\textbf{collective learning $\to$ decision making}.
\end{center}
\vskip.1cm
This dynamic will also serve as the general framework for the next two subsections of this chapter.

The following statement by Simon links the path from \textit{learning} to the process of \textit{decision} to the gain and loss of rationality within the artificial world.
\vskip.1cm
 \begin{quote}
 \textit{The human being striving for rationality and restricted within the limits of his knowledge has developed some working procedures that partially overcome these difficulties. These procedures consist in assuming that he can isolate from the rest of the world a closed system containing a limited number of variables and a limited range of consequences.}
\end{quote}
\vskip.1cm

For large systems in a collective dynamics, we can also refer to the following definition of \textit{collective learning}:

\vskip.1cm
 \begin{quote}
 \textit{Collective learning is generally defined as a social process of cumulative knowledge, based on a set of shared rules and procedures which allow individuals to coordinate their actions in search for problem solutions. This interpretation implies a conscious attitude of local agents, which find in a cooperative behavior a useful and efficient strategy to share a common knowledge.
}
\end{quote}
\vskip.1cm

An active particle approach to collective learning has been proposed in~\cite{[BDG16],[BDG16A],[BGO17]}, with further developments related to social dynamics, see~\cite{[BD19]}. The authors of these papers start from~\cite{[PIA76]} and move to~\cite{[BAN89],[SP98]}, also considering spatial dynamics and interdisciplinary applications, see~\cite{[CAP99]} and~\cite{[BA01],[Holland88],[LW98]}.
The modeling approach has been developed as follows:
\begin{enumerate}
\vskip.2cm \item \textit{Perception}: The dynamic of learning phenomena develops through \textit{Perception}, where  each individual possesses a perceptual domain within which the presence of other individuals is perceived with a different intensities. Interactions involve not only individual-based encounters, but also a process, in which each individual also learns from the collectivity viewed as a whole, which can be represented by low-order moments.

\vskip.2cm \item\textit{Collective Learning}: Interactions induce a learning process in which individuals modify (increase) their level of knowledge. Learning can then be interpreted as a \textit{cumulative dynamic}, since it accumulates over time. This dynamic may be somewhat related to memory dynamics, see~\cite{[GS96]}.

\vskip.2cm \item\textit{Decision Making}: Learning is followed by decisions. Examples were given in Section 3 concerning the dynamics of crowds, where pedestrians, viewed as a-particles, learn about the distribution of the other a-particles. The same approach can be applied to the movement of cells, examples are given in~\cite{[CKS21],[SSW21],[ZS22]}. The study of social systems and economics can focus on all the studies of social dynamics and economic systems discussed in Section 4.

\end{enumerate}

\vskip.2cm \noindent  \textbf{Remark 6.1.} \textit{Further developments could be addressed to consider interactions at different scales. Indeed, low order moments may not provide sufficient information in some cases, while learning dynamics should consider the full distribution. Then the interaction would have to consider a distance between interacting entities that includes not only the state of the particles but also the respective distribution functions. An interesting contribution to this topic is given in~\cite{[Toscani25]}. Then appropriate metrics should be defined, see~\cite{[AGS05]}.}

\subsection{Reasonings on artificial intelligence methods}\label{Sec:6.3}

In this subsection, some considerations are proposed to understand how the mathematical theory reviewed in the previous section can contribute to the problems raised by the development of the design of specific tools of Artificial Intelligence, in particular what is called Scientific Machine Learning, or SciML for short. In fact, this is a new frontier of applied mathematics that has already attracted the interest of various scientists, as can be seen from the articles published in the SIAM Journal, see~\cite{[LIN24],[NAN24]} and, of course, in the specialized literature on the subject, see~\cite{[BLH21],[CUN24],[TDC24]}, see also~\cite{[Crutchfield],[MK23]}.

An excellent review in which the interactions between numerical modeling and physics are reviewed and critically analyzed has recently been proposed in~\cite{[Alfio25]}. The theory is related to the dynamics of the heart as a specific application. This paper proposes the following definition of SciML. In fact, it is a concise but well focused description that contributes to address the research activity in the field.
\vskip.1cm
\begin{quote}
\textit{Scientific Machine Learning is an interdisciplinary field empowered by the synergy of physics-based computational models with machine-learning algorithms for scientific and engineering applications.}
\end{quote}
\vskip.1cm

The ``Nature Machine Learning''  editorial article, see~\cite{[Editorial]}, is in line with the above definition and critically analyzes the perspectives and challenges of the complex interaction between physics and PDEs, see~\cite{[Editorial]}. There is an important literature on this topic, see~\cite{[BNK23],[GKK25],[KLMQ25]}. See also~\cite{[CS22]}, for the complex interaction between decision making and neuroscience.

We  start from the perspectives introduced in~\cite{[BDL24]}, focusing on the objective of deriving a differential system capable of describing the dynamics by which the brain operates in decision making. This challenging goal is presented by Yann LeCun as a key component of the whole process of artificial intelligence, see~\cite{[CUN24]}, and we refer again to the aforementioned \textit{collective learning $\to$ decision making} dynamics.

The goal is the mathematical description of the highly complex behavior of our mind, taking into account rational and irrational behavior, as well as the complex choice between two (or more) competing strategies, see~\cite{[WC24],[Zhang2013]}, involving different actions and~\cite{[XIA15]}. In short, we are looking for a \textit{differential system capable of describing the way living things think}. Interesting clues are offered in the special issue~\cite{[Schoeller]} devoted to understanding the complex physics of the mind. Studies on this topic look forward to a systems approach to the brain, see~\cite{[DeVico],[RAF19],[Tozzi]}.

This statement, loosely interpreted, is simply naive. In fact, the complexity of the brain is such that mathematics is still far from having the tools to contribute effectively to this mythical goal. Accordingly, let us first limit our considerations to a specific dynamics, and then we can focus on some well-defined aspects of learning and decision making. We will then focus on the collective dynamics of active particles and the mathematical tools of active particle kinetic theory. The previous sections have shown that specific models can be derived to describe different types of collective dynamics. For example, in crowd dynamics it is possible to study contagion dynamics in epidemics, see~\cite{[ABKT23],[BL24]}, or the role of leaders in supporting the search for safe evacuation dynamics, see~\cite{[L23],[LZ22]}.

For both crowd and swarm dynamics, coordination leads to what has been called \textit{swarm intelligence}, see~\cite{[BHLY24],[KKHP]}. In fact, the above interaction rules coordinate the systems and lead to their collective motion, see~\cite{[BHLY24]}.
In this kind of applications the concept of \textit{swarm intelligence} appears as the collective ability of self-organization, see~\cite{[BW1989]}.
This definition can be considered as a first step to define the concept of artificial intelligence in behaviorally heterogeneous systems, where a description of the dynamics by a differential system is a necessary step to complement such a concept. Further interpretations are discussed in~\cite{[BHLY24]} with a focus on the mathematical theory of behavioral swarms, more generally it applies to all swarm models including those with internal variables, see~\cite{[HKKS19],[HKR18]}.

We are looking for a differential structure to describe the dynamics by which the human brain develops the learning decisions mentioned above. The dynamics (and the differential system) are different for each case study considered. Therefore, we first develop some general considerations and then apply them specifically to crowd dynamics.

The first step of this approach focuses on modeling collective learning, which is the precursor of many, perhaps all, collective dynamics in social dynamics and economics, but also in dynamics such as crowds and swarms. We refer to~\cite{[BDG16]}, which was developed to describe the collective dynamics of heterogeneous populations in which rational and irrational behavior coexist, see~\cite{[BAN89],[BGO17]}. Furthermore, we are guided by LeCun's concept that the amount of learning through interaction is statistically much higher than learning through individual study, see~\cite{[CUN24]}. Collective learning is commonly defined as
\vskip.1cm
\begin{quote}
\textit{Social process of cumulative knowledge based on a set of shared rules and procedures that allow individuals to coordinate their actions in the search for solutions to problems.}
\end{quote}
\vskip.1cm

The perspectives proposed in~\cite{[BDL24]} focus on the idea of designing \textit{training-predictive platforms} motivated by crisis managers for safety problems in collective dynamics of interest for our society. The rationale for their development is summarized as follows:
\begin{enumerate}
\vskip.2cm \item Selection of a specific class of  multi-agent dynamical systems and identification of the \textit{variables} that describe the dynamics and are related to learning and decision making.

\vskip.2cm \item Derivation of a \textit{differential system describing the dynamics of learning}. In fact, individuals learn at two different scales. In particular, individual learning is developed at the microscopic scale. On the other hand, individuals also learn from the global properties of the system at the macroscopic scale.

\vskip.2cm \item  Derivation of a \textit{differential system} that describes the dynamics of the decision. This dynamic is interactive with the learning that provides the inputs necessary for a fast decision. The design should include the computational codes and all the input-output devices.

\vskip.2cm \item  The differential system can be inserted in a training-prediction platform that is suitable for the description of the dynamics in different geometries and in heterogeneous emotional states.

\end{enumerate}

This topic has been developed in~\cite{[BDL24]} with a focus on safety problems in crowd dynamics. A challenging perspective is to derive a general theory that can be applied to a wide variety of collective motions of interacting living entities.

\section{Complexity towards a mathematical theory}\label{Sec:7}

This final section explores future directions for the mathematical theories discussed and analyzed earlier. An important source of ideas is the bibliography, which focuses on complexity and the collective behavior of living systems. Since our review mainly examines collective dynamics, one way to understand complexity is that the behavior of a few individual entities does not directly explain how a group behaves as a whole. The main challenge is to find mathematical equations that can bridge this gap. A good summary of the topic, including its historical development, can be found in \cite{[Krakauer]}, as well as in~\cite{[MM23],[Prig-Sten18]}, see also~\cite{[MCM00]}.

An important contribution is Giorgio Parisi's Nobel Lecture~\cite{[Parisi2023]}, which links these ideas to complex systems (from spin glasses to biological systems) that share a common mathematical structure in which microscopic disorder and nonlinear interactions give rise to emergent behavior. This approach provides a framework for interpreting the collective dynamics analyzed in this review, in particular the role of heterogeneity in learning processes.

An interesting project on complexity, still in progress, is being developed by the Santa Fe Institute. The research activity is brought to the attention of readers at the following address:
\vskip.1cm
\begin{center}
https://www.foundationalpapersincomplexityscience.org/tables-of-contents
\end{center}
\vskip.1cm

The study of collective psychology helps to understand the key features of complexity. Indeed, the psychology of large aggregations of individuals is different from that of a few individuals. A pioneering work on this topic was published in the eighteenth century~\cite{[LB1895],[LB19]}, but it was criticized because of the author's ideological position. On the other hand, we have learned that collective psychological behavior can also be manipulated through social media. The effects range from the dynamics of open markets~\cite{[CPT05],[MMT19],[Olson]} to political frameworks~\cite{[AR06],[DKLM17]}.

First, we focus on some specific features of the method. Then we report on some perspectives related to these features as well as to theories of complexity. Finally, we examine the theory as a whole to identify some possible perspectives. The key features mentioned above correspond to the following steps in the development of the theory.

\begin{itemize}

\vskip.2cm  \item The approach starts with the search for a mathematical structure that is suitable to capture the complexity of living systems without changing the general concept. The following specific topics have been identified and described in Section 2: ability to express strategy; nonlinear interactions; heterogeneity; learning ability; Darwinian mutations and selection. Others may be added or substituted, but the reasoning remains the same.

\vskip.2cm  \item General structures have been derived. These structures go far beyond the mathematical theories of statistical physics. Such structures are characterized by the fact that the dependent variable is a distribution function over the microscopic state of the interacting living entities. Such a state includes, in addition to mechanical variables, a behavioral variable called ``activity''. Some elementary structures have been proposed to model interactions, such as interaction rates, state transitions (progression/regression), and transitions between groups of aggregated entities called functional subsystems.

\vskip.2cm  \item The modeling of interactions should refer to each specific system and should be modeled under the assumption that the dynamics occur in two sequential steps. Indeed, individual entities first start learning from other individuals in their sensitivity domain, and then a decision is made with respect to the specific behavioral dynamics of the system under consideration.

\vskip.2cm  \item The mathematical tools reviewed in our paper have always been motivated by the search of a mathematics for living systems. The concept that has been developed is that the strategy should have the ability to preserve the complexity features of living systems. These where described in Subsection 2.2. These features are also an important reference for methods such as~\cite{[NOW06]}, which differ with respect to those treated in our paper, but work for the mythical goal of describing the dynamics of living matter.

\vskip.2cm  \item The presentation has also included a brief review of methods that have been independently developed in parallel. These methods have generated a rich literature and motivated the search for a unified mathematical theory which, if effectively achieved, would make a decisive contribution to the mythical objective mentioned above.
\end{itemize}

Before moving on to perspectives, we should try to provide an answer to the key question posed in subsection 2.3, i.e. \textit{To what extent can the mathematical tools to be derived according to the above general framework be considered a mathematical theory?}. We can argue that the differential system could be studied within the framework of the philosophy of mathematics. An excellent reference is the book~\cite{[CC22]}. Therefore, we can only say that we have developed a theoretical approach to derive a general differential system that can be related to various complex systems. This quest was also motivated by pages 469-470 of Simon's seminal paper, see~\cite{[Simon1962]}.

A forward look at research perspectives, selected according to the research expertise of the authors, it could even be the bias of the authors, are brought in the notes below. A brief description is given for each, while a concluding comment shows how these key issues are all interrelated.

\vskip.2cm  \noindent \textbf{Learning dynamics:} All the dynamics we have discussed in this paper start with an initial state in which a-particles learn from the others. This process is consistent with Krakauer's information-theoretic perspective, where learning in complex systems can be modeled as information flow across hierarchical levels, with individual strategies emerging from thermodynamic and cognitive constraint~\cite{[KBOFA20]}. The subsequent dynamics incorporate previously acquired information, a phenomenon observed in traditional teaching-learning systems~\cite{[BA01],[BDG16A]} and social interactions~\cite{[Holland88],[LW98],[SP98]}. The Editorial in Nature reports about the current strategies, and conceptual difficulties,
of the development of Scientific Machine Learning, see~\cite{[Editorial]}.

Recent advances show that these dynamics are increasingly relevant to technological processes~\cite{[BDKV20],[DOSI2023],[DPV17]} and deep learning~\cite{[BLH21]}. As Krakauer's work shows, the heterogeneity of learning agents fundamentally shapes collective outcomes - a key insight for modeling multiscale systems where individual microscopic learning interacts with macroscopic emergent properties.

Social interactions involve interdisciplinary dialogue, such as medicine, epidemiology, politics, and economics, which is even more critical in times of pandemics. This was the message delivered in~\cite{[BBC20]} and subsequently considered by various authors, as examples see~\cite{[ADKV],[B3EPT],[CAP99]}.

This brief overview of the bibliography is only a collection of examples and does not claim to be complete. However, it can be stated that the kinetic theoretical approach developed in~\cite{[BDG16]} was a visionary anticipation of the current mathematical studies of the various types of collective dynamics. As we have seen in the previous sections, learning involves a rich variety of interactions that evolve across scales and different scientific environments. We believe that future research should move in this direction.

\vskip.2cm  \noindent \textbf{Multiple strategies and the utility functions:} The heterogeneity of living systems appears not only in the learning dynamics, but also in the decision dynamics, where different strategies can appear motivated by the \textit{utility function}. Let us consider that the last chapter of the book~\cite{[BL12]} emphasizes that the study of social systems should go beyond the dynamics of consensus, while the review~\cite{[ABG16]} proposes that consensus and dissent coexist in a heterogeneous system in which individuals, in which individuals chose one of the two strategies based on their own state and, if the case, under the influence of external actions.

The so-called Parrondo's paradox is critically analyzed in a recent paper that focuses on the role of two strategies that lose when used individually, but can win when used in an appropriate combination, see~\cite{[WC24]} and the comment~\cite{[NO25]}.

If a the mathematical approach is developed by differential systems derived according to Simon's philosophy, see~\cite{[BE24]}, the mathematical interpretation is that the utility function is not constant in time, but it is heterogeneously distributed over the population and evolves in time thin the  \textit{artificial world} made of \textit{artifacts} created by the interacting individuals. The utility function also depends depends on the ability of each individual entity (active particle) to look forward and make what is the best individual \textit{payoff} also accounting for a long time vision ability.

\vskip.2cm  \noindent \textbf{Multiscale dynamics towards a unified mathematical theory:} The need of a multiscale view has been mentioned in the previous sections motivated by Hilbert's 6th problem~\cite{[Hilbert]}. The main problem is the derivation of models at the macroscopic scale from the underlying description at the microscopic scale. Some literature has been cited in Subsection 3.3 to show how this problem has been dealt with in some specific case studies.

This topic can contribute to the search for a unified mathematical theory that is capable of taking into account not only the contents of Sections 3 and 4, but also some contributions briefly reviewed in Section 5. Furthermore, the considerations in subsection 6.3, referring to mathematical tools for artificial intelligence, indicate that below the microscopic scale, dynamics needs at least one further scale corresponding to the specific class of systems under consideration. Examples in biology point to the molecular scale of genes, while the search for an artificial brain points to the complexity of neural networks. This lower scale plays an important role in both learning and decision dynamics.

\section*{Acknowledgments}

\noindent Nicola Bellomo acknowledges the partial support by the State Research Agency of the Spanish Ministry of Science and FEDER-EU, project PID2022-137228OB-I00 (MICIU/AEI /10.13039/501100011033); by Modeling Nature Research Unit, Grant QUAL21-011 funded by Consejer\'ia de Universidad, Investigaci\'on e Innovaci\'on (Junta de Andaluc\'ia).

\vskip.2cm \noindent Jie Liao  acknowledges the partial support by The National Natural Science Foundation of China 12471211.

\vskip.2cm \noindent Nicola Bellomo, Diletta Burini, and Jie Liao, would like to express their gratitude to their colleagues Massimo Egidi and Pietro Terna, for their precious suggestions in the field of behavioral economics, complexity and Herbert A. Simon's philosophy.

%
%
%

\end{document}